\documentclass[12pt]{article}
\usepackage{amsmath}
\usepackage{amssymb}
\usepackage{amscd}

\newtheorem{Df}{Definition}[section]
\newtheorem{Te}[Df]{Theorem}
\newtheorem{Po}[Df]{Proposition}
\newtheorem{Cr}[Df]{Corollary}
\newtheorem{Lm}[Df]{Lemma}
\newtheorem{Ca}[Df]{Claim}
\newtheorem{Cn}[Df]{Conjecture}

\newcommand{\Bdf}{\begin{Df}}
\newcommand{\Edf}{\end{Df}}
\newcommand{\Bte}{\begin{Te}}
\newcommand{\Ete}{\end{Te}}
\newcommand{\Bpo}{\begin{Po}}
\newcommand{\Epo}{\end{Po}}
\newcommand{\Bcr}{\begin{Cr}}
\newcommand{\Ecr}{\end{Cr}}
\newcommand{\Blm}{\begin{Lm}}
\newcommand{\Elm}{\end{Lm}}
\newcommand{\Bca}{\begin{Ca}}
\newcommand{\Eca}{\end{Ca}}
\newcommand{\Bcn}{\begin{Cn}}
\newcommand{\Ecn}{\end{Cn}}
\newcommand{\Bdm}{{\it Proof.}\ }
\newcommand{\Edm}{\rule{2mm}{2mm}}
\newcommand{\Rm}{{\it Remark \arabic{section}.\arabic{Df} \ }}
\newcommand{\Ea}{{\it Example \arabic{section}.\arabic{Df} \ }}

\begin{document}

\title{\bf{Koszul and Gorenstein properties for homogeneous algebras}}
\author{\bf{Roland Berger} \textrm{and} \bf{Nicolas Marconnet}\\ \\ \emph{LARAL},
Facult\'{e} des Sciences et Techniques\\ 23, Rue P. Michelon,  42023 Saint-Etienne Cedex,
France \\
\emph{E-mail}: Roland.Berger@univ-st-etienne.fr \\ Nicolas.Marconnet@univ-st-etienne.fr}
\date{}
\maketitle
\begin{abstract}
Koszul property was generalized to homogeneous algebras of degree $N>2$
in~\cite{rb:nonquad}, and related to $N$-complexes in~\cite{bdvw:homog}. We show that if
the $N$-homogeneous algebra $A$ is generalized Koszul, AS-Gorenstein and of finite
global dimension, then one can apply the Van den Bergh duality theorem~\cite{vdb:dual} to
$A$, i.e., there is a Poincar\'e duality between Hochschild homology and cohomology of $A$, as for
$N=2$. 

\end{abstract}

\section{Introduction}
The appearance of noncommutative projective geometry has given a new impulse 
to the study of (noncommutative) graded algebras~\cite{st:npg}. Some interesting classes of graded
algebras have been investigated, for example, the class of the AS-regular
algebras~\cite{as:regular}. Those algebras are often defined as graded algebras of
finite global dimension and having AS-Gorenstein property (the polynomial
growth imposed by Artin and Schelter is often removed and in fact, it is not necessary for our
purpose). Generalized Koszulity is another property which is natural to assume for getting
noncommutative analogues of polynomial algebras~\cite{rb:nonquad, bdvw:homog}. This article examines
the class of AS-regular algebras which are generalized Koszul. Our main result is that the
Poincar\'e duality holds in this class. Besides the quadratic case
already known by Michel Van den
Bergh~\cite{vdb:dual}, the first known  Poincar\'e duality in a nonquadratic homogeneous situation
has been obtained by the second author~\cite{nic:cubic} for the
generic cubic AS-regular algebras (with polynomial growth) of global dimension three and of type A. 

The graded algebras of our class are $N$-homogeneous where $N$ is an integer $\geq 2$, that is,
they are generated in degree one and all their relations are homogeneous of the same degree $N$.
Koszul property is well-known for quadratic algebras (i.e., for $N=2$)~\cite{pri:kos}, and it was
recently generalized by the first author~\cite{rb:nonquad} for
$N>2$. This generalization follows along the definition given by Beilinson, Ginzburg,
Soergel~\cite{bgso:kdp} for $N=2$ : the minimal projective resolution of the
trivial module has to be pure in the lowest degrees. Here the minimal projective resolution makes
sense in the category of graded modules (with morphisms of degree zero), so the definition is very
natural as far as graded algebras are concerned. On the other hand, following Manin's monograph on
quadratic algebras~\cite{man:kalg}, Michel Dubois-Violette, Marc Wambst and the first author have
given an alternative set-up of generalized Koszul property based on
$N$-complexes~\cite{bdvw:homog}. Actually, if one follows the categorical point of view of Manin,
the Koszul complex of the quadratic case becomes naturally an $N$-complex in the $N$-case.

What is an $N$-complex ? It is just a complex for which the condition $d^{2}=0$ on the
differential $d$ is replaced  by $d^{N}=0$. Thus a usual complex is a 2-complex. There exist some
reserved feelings concerning the interest of $N$-complexes when $N>2$. It lies on the fact that
$N$-homologies in various classical contexts (simplicial, Hochschild, cyclic) do not bring new
things~\cite{dv:genhom, kw:Nhom, w:Ncycl}. It perhaps comes from the fact that acyclicity of $N$-complexes is
too strong. For example, acyclicity of the Koszul $N$-complex implies that the algebra is trivial.
As shown in~\cite{bdvw:homog}, there is one (and only one) adequate contraction of the Koszul
$N$-complex whose acyclicity provides the (highly non-trivial !) Koszul property. Similarly,
AS-Gorenstein property of generalized Koszul algebras coincides with acyclicity of one contraction
of the dual Koszul $N$-complex.  

So the Koszul 2-complexes of $N$-homogeneous algebras are adequate contractions of
the Koszul $N$-complexes. Of course, the Koszul 2-complexes can be defined without
$N$-complexes, and it was doing so in~\cite{rb:nonquad}. Why is it interesting to have
$N$-complexes behind the Koszul 2-complexes ? A striking answer is the
following. If we want to define the \emph{bimodule} Koszul 2-complex of an $N$-homogeneous algebra for $N>2$,
the knowledge of the left and right Koszul 2-complexes is not sufficient. Actually, it is
impossible to avoid $N$-complexes in defining the differential of the bimodule Koszul
2-complex~\cite{rb:nonquad2}. The precise reason is that the left and right $N$-differentials
commute, but not the left and right 2-differentials.

The usefulness of $N$-complexes can also be demonstrated from
constructions of Theoretical Physics (see~\cite{bdvw:homog,dvp:plac} and references therein). It is worth
noticing here that Yang-Mills algebras studied recently by Connes and
Dubois-Violette~\cite{cdv:yang} provide new examples of cubic Koszul AS-regular algebras.

Let us present the plan and some results of this paper. In Section 2, generalized Koszulity is
introduced as in~\cite{rb:nonquad}, including several points which were omitted
in~\cite{rb:nonquad}. Koszulity of Manin products is discussed.

In Section 3, the link with
$N$-complexes is explicited. By shifting the Koszul $N$-complex, the Yoneda algebra of a generalized
Koszul algebra is computed, with particular attention to the signs in the Yoneda product. 

In Section 4, the bimodule
Koszul 2-complex is drawn from the bimodule Koszul
$N$-complex, and is used to get a bimodule characterization of Koszulity. A consequence is that
Hochschild dimension and global dimension of generalized Koszul algebras are the same. 

In Section 5, AS-Gorenstein property is recalled. If $A$ is generalized Koszul over the ground
field $k$, with finite global dimension $D$ such that Ext$_{A}^{D}(k,A)=k$, then an explicit
$N$-complexes morphism
$\Phi$ is constructed from the dual Koszul $N$-complex to the Koszul $N$-complex. From $\Phi$, we
deduce a convenient criterion for $A$ to be AS-Gorenstein (Theorem 5.4). This criterion
has the following consequences.
\Bte
For $n\geq N\geq 2$, consider the $N$-homogeneous algebra $A$ over $k$ of characteristic 0, with
generators
$x_{1},\ldots ,x_{n}$. The relations of $A$ are the antisymmetrizers of degree $N$, i.e., are the
following
$$\sum_{\sigma \in \mathbf{S}_{N}} \mathrm{sgn} (\sigma)\ x_{i_{\sigma(1)}}\ldots
x_{i_{\sigma(N)}}=0\, ,\ 1\leq i_{1}<\cdots <i_{N} \leq n\,,$$
where $\mathbf{S}_{N}$ denotes the permutation group of $\{1, \ldots ,N\}$. We
know~\cite{rb:nonquad} that $A$ is generalized Koszul, of finite global dimension. Then $A$ is
AS-Gorenstein if and only if either $N=2$ or \emph{($N>2$ and $n= Nq + 1$ for some integer
$q\geq 1$)}.
\Ete  
\Bte
Let $A$ be $N$-homogeneous, generalized Koszul, of finite global dimension. Denote by $E(A)$ the
Yoneda algebra of $A$. Then $A$ is AS-Gorenstein if and only if the
finite dimensional algebra $E(A)$ is Frobenius.
\Ete

The latter result is due to Smith~\cite{smi:some} when $N=2$. Our proof follows along the same
lines, but the explicit $\Phi$ is not in~\cite{smi:some}. The situation $N>2$ is a bit more involved
because the jump map defining the adequate contraction of the Koszul $N$-complex is not additive.
Furthermore we have to pay some attention to the signs in the Yoneda product. Actually, Theorem 1.2 is a straightforward consequence of a general theorem by Lu, Palmieri, Wu and
Zhang~\cite{lpwz:ainf}, which asserts that if $A$ is any connected graded algebra, then $A$ is
AS-regular if and only if the Yoneda algebra of $A$ is Frobenius. Their
proof uses $A_{\infty}$-algebras. Under generalized Koszulity assumption, $A_{\infty}$-algebras
are not necessary. But the $N$-complexes machinery is required if we do not want to use the
general theorem. We thank Bernhard Keller to have pointed out to us
the article~\cite{lpwz:ainf}.

In Section 6, we prove the assumptions of the Van den Bergh duality theorem~\cite{vdb:dual} for our
class of algebras. An essential ingredient is the adaptation to the $N$-case of a
computation used by Van den Bergh in the situation of the quadratic Koszul noetherian
AS-regular algebras (Theorem 9.2 in~\cite{vdb:exist}). Note that our algebras may be non
noetherian (those of Theorem 1.1 are not noetherian if $N>2$), so that the Yekutieli balanced
dualizing complex~\cite{yek:dualizing} cannot be used. Then we get our
main result ($HH$ denotes Hochschild homology or cohomology) :
\Bte
Let $A$ be an $N$-homogeneous algebra which is generalized Koszul, AS-Gorenstein and with finite
global dimension $D$. For any bimodule $M$, we have
$$HH^{i}(A,M)\cong HH_{D-i}(A,\,_{\varepsilon^{D+1}\phi}M).$$
\Ete

For $N=2$, this is Van den Bergh's Proposition 2~\cite{vdb:dual}. The meaning of $\varepsilon$ is
the same as for
$N=2$. Let $\nu$ be the automorphism of the Frobenius algebra
$E(A)$ (Theorem 1.2 above) such that the bimodule $E(A)^{\ast}$ is canonically $E(A)_{\nu}$
(usual action on the left and action twisted by $\nu$ on the right). The automorphism $\phi$ of $A$
is such that the homogeneous component
$\phi_{1}$ is the transposed linear map of $\nu_{1}$, as for $N=2$. It
is important to note that, if $N>2$, $E(A)$ is
not generated in degree 1, and $\nu$ is not determined by $\nu_{1}$ (but $E(A)$
is generated in degree 1 and 2, and $\nu$ is determined by $\nu_{1}$
and $\nu_{2}$). The end of Section 6 is devoted to give direct
proofs of Theorem 1.3 for $M=A$, and specially to recover the proof used by the second author for the generic cubic
AS-regular algebras (with polynomial growth) of global dimension three and of type
A~\cite{nic:cubic}. This
proof works well for any $A$ such that there is no twist on $M=A$,
i.e., such that $\phi=\varepsilon^{D+1}$. On the other hand,
as for $N=2$, the automorphisms $\nu$ and $\phi$ are interpreted in
terms of the Artin-Schelter matrix $Q$~\cite{as:regular} when $A$ is a
cubic AS-regular algebra (with polynomial growth) of
global dimension three. 
\setcounter{equation}{0}

\section{Generalized Koszul algebras}

Some material about generalized Koszul property can be found
in~\cite{rb:nonquad, rb:nonquad2, bdvw:homog, cdv:yang, dvp:plac,
  nic:cubic}. It is completed here on several particular
points. Throughout this paper, $k$ is a field and $V$ is a finite-dimensional $k$-vector space
which is considered as graded, with its gradation concentrated in degree 1. We fix a natural number
$N\geq 2$ and a graded subspace $R$ of the tensor algebra Tens$(V)$, which is concentrated in degree
$N$. In other words $R$ is a subspace of $V^{\otimes N}$. The two-sided ideal $I(R)=I$ generated by
$R$ in Tens$(V)$ is graded by the subspaces $I_{n}$ given by $I_{n}=0$, $0\leq n
\leq N-1$, and 
$$I_{n}=\sum_{i+j+N=n}V^{\otimes i}\otimes R\otimes V^{\otimes j}\,,\ n\geq N.$$
The algebra $A$=Tens$(V)/I$ is called an \emph{N-homogeneous algebra} on $V$, with $R$ as space of
relations. The algebra $A$ is graded by the subspaces $A_{n}=V^{\otimes n}/I_{n}$, and generated by
$V$ (hence in degree 1). Clearly, $A_{n}=V^{\otimes n}$ for $0\leq n \leq N-1$. Introduce the
jump map $n : \mathbf{N}\rightarrow \mathbf{N}$ by 
$$n(2q)=qN, \ n(2q+1)= qN+1,\ q\ \mathrm{integer} \geq 0.$$ 
The jump map $n$ is additive if and only if $N=2$. If $N>2$, $n(i+j)=n(i)+n(j)$ holds unless if $i$
and $j$ are both odd. We begin to generalize a result of~\cite{bgso:kdp}.
\Bpo
Let $A$ be an $N$-homogeneous algebra. The graded vector space $\mathrm{Tor}_{i}^{A}(k,k)$ lives in
degrees
$\geq n(i)$ for $i\geq 0$.
\Epo
\Bdm
We prove the result under the more general assumption that $R$ lives in degrees $\geq N$. Let
$P\rightarrow k$ be a minimal projective resolution of the trivial module $k$ in the category
$A$-grMod of the graded left $A$-modules (with morphisms of degree 0). We know that
$P_{i}=A\otimes E_{i}$, where $E_{i}$ is a graded vector subspace of $P_{i-1}$, and the
differential $d_{i}: P_{i} \rightarrow P_{i-1}$ is defined by this inclusion. Since the
graded vector space Tor$_{i}^{A}(k,k)$ is isomorphic to
$k\otimes_{A} P_{i}=E_{i}$, it suffices to prove that $E_{i}$ lives in degree $\geq n(i)$
for $i\geq 0$. Proceed by induction on $i$. We can choose $E_{0}=k$,
$E_{1}=V$, $E_{2}=R$ with obvious inclusions, so that the property
holds for $i=$ 0, 1, 2. Assume
$i\geq 3$ and the property true for $i-1$. First assume $i=2q+1$. Since $P_{i-1}$
lives in degrees $\geq qN$, $E_{i}$ lives in degrees $\geq qN$. But $P_{i-1,qN}=E_{i-1,qN}$, so
that $\ker(d_{i-1})$ vanishes in degree $qN$. Hence, $E_{i,qN}=0$ and $E_{i}$ lives in degrees $\geq
qN+1$. Assume now $i=2q$. Since $P_{i-1}$ lives in degrees $\geq (q-1)N+1$, $E_{i}$ lives in
degrees $\geq (q-1)N+1$. Fix $u$ with $1\leq u \leq N-1$. We have
$$P_{i-1,(q-1)N+u}= \bigoplus_{1\leq v\leq u} A_{u-v}\otimes E_{i-1,(q-1)N+v}\,.$$
For proving that $E_{i}$ lives in degrees $\geq qN$, it suffices to prove that $d_{i-1}$ is
injective on each $A_{u-v}\otimes E_{i-1,(q-1)N+v}$. But $d_{i-1}$ restricted to
$E_{i-1,(q-1)N+v}$ is the inclusion into the direct sum
$$\bigoplus_{1\leq w\leq v} A_{v-w}\otimes E_{i-2,(q-1)N+w}\,.$$
Therefore, $d_{i-1}$ sends $A_{u-v}\otimes E_{i-1,(q-1)N+v}$ into the direct sum
$$\bigoplus_{1\leq w\leq v} A_{u-w}\otimes E_{i-2,(q-1)N+w}\,.$$
Since $u-w<N$, $A_{u-w}=V^{\otimes (u-w)}$ which is not a quotient, hence the result.
\Edm
\Bdf
An $N$-homogeneous algebra $A$ is said to be (generalized)
Koszul if for each $i\geq 3$, the graded vector space \emph{Tor}$_{i}^{A}(k,k)$ is
concentrated in degree $n(i)$ (or, equivalently, \emph{\underline{Ext}}$_{A}^{i}(k,k)$ is concentrated in degree
$-n(i)$).
\Edf

When $N=2$, it is exactly Priddy's definition~\cite{pri:kos}. Here, \underline{Ext} denotes the derived
functor of the functor \underline{Hom} of $A$-grMod. Let us denote by hom the functor Hom
of $A$-grMod (in which the morphisms are of degree 0), and ext its
derived functor. For any objects $M$ and $M'$ of $A$-grMod, we have
\underline{Ext}$_{A}^{i,j}(M,M')=$ ext$_{A}^{i}(M,M'(j))$. The following
proposition generalizes a result of~\cite{bgso:kdp}. Recall that the gradation of the shift $M(l)$
of a graded 
$A$-module $M$ is defined by $M(l)_{n}=M_{n+l}$.
\Bpo
Let $A$ be a connected graded $k$-algebra. Fix $N\geq 2$. The following are equivalent.

(i) $A$ is $N$-homogeneous and Koszul.

(ii) \emph{ext}$_{A}^{i}(k,k(-n))=0$ whenever $n\neq n(i)$.

(iii) For any objects $M$ and $N$ of $A$-grMod, respectively concentrated in degrees $m$ and $n$,
\emph{ext}$_{A}^{i}(M,N)=0$ whenever $n\neq m+n(i)$.
\Epo
\Bdm
Clearly, \underline{Ext}$_{A}^{i,-n}(k,k)=$ ext$_{A}^{i}(k,k(-n))$, hence (i)$\Leftrightarrow$(ii).
(iii)$\Rightarrow$(i) is obvious. Assume (i). Let $P$ be a minimal projective resolution of $k$,
$P_{i}=A\otimes E_{i}$ with $E_{i}$ concentrated in degree $n(i)$. As
ext$_{A}^{i}(M,N)=$ ext$_{A}^{i}(M(-m),N(-m))$, we assume
$m=0$. As $M$ is a direct sum of copies
of $k$, we assume $M=k$. Since hom$_{A}(P_{j},N)=$ hom$_{k}(E_{j},N)$ vanishes whenever $n\neq
n(j)$, each component of the complex hom$_{A}(P,N)$ vanishes unless perhaps one. Therefore,
ext$_{A}^{i}(k,N)=$ hom$_{k}(E_{i},N)$, which vanishes whenever $n\neq n(i)$.
\Edm

As for $N=2$, there exists a remarkable complex whose acyclicity is equivalent to Koszulity of
$A$. For $n\geq 0$, introduce the subspace of $V^{\otimes n}$ by
$$W_{n}=\bigcap_{i+j+N=n}V^{\otimes i}\otimes R\otimes V^{\otimes j}.$$
If $n<N$, $W_{n}=V^{\otimes n}$. As defined in~\cite{rb:nonquad}, the \emph{Koszul 2-complex} $K$
of $A$ is the following complex of $A$-grMod
\begin{equation}
\cdots \longrightarrow K_{i} \stackrel{\delta_{i}}{\longrightarrow} K_{i-1} \longrightarrow \cdots
\longrightarrow K_{1} \stackrel{\delta_{1}}{\longrightarrow} K_{0} \longrightarrow 0\,,
\end{equation}
in which $K_{i}=A\otimes W_{n(i)}$, and the differential $\delta_{i}$ is defined by the inclusion
of $W_{n(i)}$ in $K_{i-1}$. The homology of $K$ is $k$ in degree 0, and 0 in degree 1.
\Bte
Let $A$ be $N$-homogeneous. Then $A$ is Koszul if and only if the
Koszul $2$-complex $K$ is exact in
degrees $\geq 2$, i.e., is a resolution (necessarily minimal projective in $A$-grMod) of $k$.
\Ete
\Bdm
If $K$ is a resolution of $k$, then Tor$_{i}^{A}(k,k)=W_{n(i)}$, which is concentrated in degree
$n(i)$, hence $A$ is Koszul. Conversely, assume $A$ Koszul. We want to prove $H_{i}(K)=0$ by
induction on $i\geq 1$. It holds for $i=1$. Let $i\geq 2$ be such that it holds for any integer
$<i$. Set $Z_{i}=\ker \delta_{i}$. By dimension
shifting (\cite{weib:homo} p.47), one has
$$\mathrm{Tor}_{i+1}^{A}(k,k)= \ker (k\otimes_{A}Z_{i} \rightarrow k\otimes_{A} K_{i}).$$
But $k\otimes_{A} Z_{i}$ lives in degree $>n(i)$, whereas $k\otimes_{A}K_{i}$ lives in degree
$n(i)$. Therefore, Tor$_{i+1}^{A}(k,k)= k\otimes_{A} Z_{i}$, and $k\otimes_{A} Z_{i}$ is
concentrated in degree $n(i+1)$. Thus the projective cover $A\otimes_k(k\otimes_{A}Z_{i})$ of $Z_{i}$ is generated in degree
$n(i+1)$, and it is the same for $Z_{i}$ itself. 

To conclude that $Z_{i}=\mathrm{im}\, \delta_{i+1}$, it suffices to prove this equality in degree
$n(i+1)$. But $K_{i,n(i+1)}$ is $V\otimes W_{n(i)}$ or $V^{\otimes (N-1)}\otimes W_{n(i)}$,
according to $i$ is even or odd. Consequently $\delta_{i+1}(W_{n(i+1)})$ is the natural
inclusion of $W_{n(i+1)}$ into $V\otimes W_{n(i)}$ or $V^{\otimes (N-1)}\otimes W_{n(i)}$. On the
other hand,
$\delta_{i}$ sends naturally $V\otimes W_{n(i)}$ or $V^{\otimes (N-1)}\otimes W_{n(i)}$ into
$A_{N}\otimes W_{n(i-1)}$. If $i$ is even, we have
$$(\ker \delta_{i})_{n(i+1)}=(V\otimes W_{n(i)})\cap(R\otimes W_{n(i-1)}),$$
which is $W_{n(i+1)}$, as expected. If $i$ is odd, we have
$$(\ker \delta_{i})_{n(i+1)}=(V^{\otimes (N-1)}\otimes W_{n(i)})\cap(R\otimes W_{n(i-1)}).$$
If $N=2$, it is $W_{n(i+1)}$, as expected. It is still $W_{n(i+1)}$ if $N>2$, because, according to
the lemma just below, Koszulity of $A$ implies the following
$$(V^{\otimes (N-1)}\otimes R)\cap (R\otimes V^{\otimes
  (N-1)})=W_{2N-1}.\ \Edm $$
\Blm
If $A$ is Koszul $N$-homogeneous, one has for $m=1,\ldots, N-1$,
\begin{equation}
(V^{\otimes m}\otimes R)\cap (R\otimes V^{\otimes m})=W_{N+m}.
\end{equation}  
\Elm
\Bdm
In the previous induction, we can take the equality $H_{i}(K)=0$ for granted when $i=2$. So
$(\ker\delta_{2})_{N+m}= V^{\otimes (m-1)}\otimes W_{N+1}$. Since
$\delta_{2}$ is the natural map from $V^{\otimes m} \otimes R$ to
$A_{N+m-1} \otimes V$, we get 
\begin{equation}
(V^{\otimes m}\otimes R)\cap (R\otimes V^{\otimes m}+\cdots + V^{\otimes (m-1)}\otimes R
\otimes V) = V^{\otimes (m-1)}\otimes W_{N+1}
\end{equation} 
which implies the inclusion
\begin{equation}
(V^{\otimes m}\otimes R)\cap (R\otimes V^{\otimes m}) \subseteq V^{\otimes (m-1)}\otimes R \otimes
V,
\end{equation} 
and this inclusion gives (2.2) for $m=2$. Assume now (2.2) for $m$. Then (2.4) for $m+1$ shows
that $(V^{\otimes (m+1)}\otimes R)\cap (R\otimes V^{\otimes (m+1)})$ is included into
$$(V^{\otimes (m+1)}\otimes R)\cap (R \otimes V^{\otimes (m+1)})\cap(V^{\otimes m}\otimes R \otimes
V).$$ Using induction hypothesis, the latter subspace is 
$$(V^{\otimes (m+1)}\otimes R)\cap (W_{N+m} \otimes V) = W_{N+m+1}.\ \Edm$$
\Bcr
Let $A$ be Koszul $N$-homogeneous. Knowing that the resolution $K$ is
minimal in $A$-grMod and setting $W_{n(\infty)}=0$, we have
$$\mathrm{gl.dim}(A) = \sup \{i\in \mathbf{N}\cup \{\infty\}\,;\ W_{n(i)}\neq 0\,\}.$$
\Ecr

Relations (2.3) for $m=2,\, \ldots,N-1$ are put together in order to define the so-called
\emph{extra condition} (actually, the relation (2.3) for $m=N-1$ implies the other ones). If the
extra condition is added to a certain infinite sequence of distributivity relations (generalized
Backelin's relations), one gets an explicit characterization for $A$ to be Koszul
(\cite{rb:nonquad}, Theorem 2.11). The extra condition is void if
$N=2$, and it is a consequence of the jump of degrees if $N>2$. Generalized Backelin's relations
are satisfied if $A$ is \emph{distributive}. Recall this definition~\cite{rb:nonquad}. Denote by
$\mathcal{L}(V^{\otimes n})$ the lattice (for inclusion) of the vector subspaces of $V^{\otimes
n}$. Infimum and supremum in this lattice are respectively intersection and sum of subspaces.
\Bdf
Let $A$ be $N$-homogeneous. One says that $A$ is distributive if for any $n\geq N$ ($n\geq N+2$
suffices), the sublattice of $\mathcal{L}(V^{\otimes n})$ generated by the subspaces $V^{\otimes
i}\otimes R\otimes V^{\otimes j}$, $i+j+N=n$, is distributive.
\Edf

Note that if $A$ is distributive, Koszulity (i.e., the extra condition) is equivalent to relations
(2.2) (or to relations (2.4)). A stronger condition (very useful in practice) than distributivity is
Bergman's confluence relatively to a basis $X$ of $V$~\cite{rb:nonquad}. For example, the algebras
in Theorem 1.1 (Section 1) are confluent for any basis of $V$. Another example of $X$-confluence is the case where the relations of $A$ are monomials in $X$. In
this case, the extra condition has a remarkable combinatorial
characterization (\cite{rb:nonquad}, Proposition 3.8). 

Let us discuss Koszulity of Manin products, which is well-known for $N=2$~\cite{bafr:kalg}. Let
$A$ and $A'$ be $N$-homogeneous algebras on $V$ and $V'$, with spaces of relations $R$ and $R'$,
respectively. Their Manin products~\cite{bdvw:homog} are the $N$-homogeneous algebras $A\circ A'$
and $A\bullet A'$,  both on $V\otimes V'$, with respective spaces of relations $R\otimes
V^{'\otimes N} +V^{\otimes N} \otimes R'$ and $R\otimes R'$. The spaces of relations are in
$(V\otimes V')^{\otimes N}$ by the canonical isomorphism $V^{\otimes
  N} \otimes V^{'\otimes N}
\cong (V\otimes V')^{\otimes N}$. Show that if $A$ and $A'$ are distributive, then $A\circ A'$
and $A\bullet A'$ are distributive. In fact (we omit tensor products to simplify), if $B$ (resp.
$B'$) is a basis of $V^{n}$ (resp. $V^{'n'}$) distributing the family $V^{i} R
V^{j}$, $i+j+N=n$ (resp. $V^{'i'} R' V^{'j'}$, $i'+j'+N=n'$), $BB'$
is a basis of $V^{n} V^{'n'}$
distributing the family  $V^{i}RV^{j}V^{'i'}R'V^{'j'}$,
$V^{i}RV^{j}V^{'n'}$,
$V^{n}V^{'i'}R'V^{'j'}$, $i+j+N=n$, $i'+j'+N=n'$.
\Bpo
Let $A$ and $A'$ be $N$-homogeneous and distributive. If $A$ and $A'$ are Koszul, then
$A\bullet A'$ is Koszul.
\Epo
\Bdm
Set $\mathcal{V}=V\otimes V'$. Consider $\mathcal{R}$ as being $R\otimes R'$ naturally viewed in
$\mathcal{V}^{\otimes N}$. Omit again tensor products in the notations. It suffices to prove the
inclusion
$$(\mathcal{V}^{N-1}\mathcal{R})\cap (\sum_{i+j=N-1,j\geq 1} \mathcal{V}^{i} \mathcal{R}
\mathcal{V}^{j} ) \subseteq \mathcal{V}^{N-2} \mathcal{R} \mathcal{V}.$$
Identifying $\mathcal{V}^{i} \mathcal{R} \mathcal{V}^{j}$ to $V^{i}RV^{j}V^{'i}R'V^{'j}$ and using
distributivity, we see that the left-hand side of the inclusion coincides with
$$\sum_{i+j=N-1,j\geq 1} (V^{N-1}R\cap V^{i}RV^{j})(V^{'N-1}R'\cap V^{'i}R'V^{'j}),$$
which is included into
$$(\sum_{i+j=N-1,j\geq 1} V^{N-1}R\cap V^{i}RV^{j})(\sum_{i'+j'=N-1,j'\geq 1}V^{'N-1}R'\cap
V^{'i'}R'V^{'j'}).$$
By Koszulity of $A$ and $A'$, the latter space is included into $V^{N-2}RVV^{'N-2}R'V'$ which is
identified to $\mathcal{V}^{N-2} \mathcal{R} \mathcal{V}$.
\Edm

For $N=2$, if $A$ is Koszul, then the dual algebra $A^{!}$ is Koszul. But it is no longer true for
$N>2$~\cite{rb:nonquad}. The same phenomenon occurs between the dual Manin products $\circ$ and
$\bullet$.

\addtocounter{Df}{1}
\Ea
Here $N=3$, $A$ (resp. $A'$) has two generators $x_{1}, x_{2}$ (resp. $x'_{1}, x'_{2}$) with the
relation $x_{1}^{3}=0$ (resp. $x_{1}^{'3}=0$). Then $A$ and $A'$ are
monomial, and the combinatorial condition for Koszulity holds for them. So $A$ and
$A'$ are Koszul. But $A\circ A'$ is not Koszul. In fact $A\circ A'$ is monomial, with set $E$ of
monomials formed of the $x_{1}^{3}m'$, $m'$ word of length 3, and of the
$mx_{1}^{'3}$, $m$ word of length 3. However the combinatorial condition for
Koszulity does not hold for $E$: 
$x_{1}^{3}x_{2}^{2}x_{2}^{'2}x_{1}^{'3}$ is such that the first and the last factor of length 3 are
in $E$, but the central factor $x_{1}^{2}x_{2}x'_{2}x_{1}^{'2}$ is not in $E$.

\setcounter{equation}{0}

\section{Koszul $N$-complexes}
Throughout this section, $A$ is an $N$-homogeneous algebra on $V$, with $R$ as space of relations.
Let us define the left Koszul $N$-complex $K_{l}(A)$ of $A$ in an elementary way (for a more
conceptual way, see~\cite{bdvw:homog}), which will be enough for us. The right Koszul
$N$-complex $K_{r}(A)$ would be defined analogously. Recall that $A^{!}$ is the $N$-homogeneous
algebra on the dual vector space $V^{\ast}$, with $R^{\perp}$ as space of relations. Introduce the
canonical element $\xi_{l}= \sum_{i} e_{i}\otimes e_{i}^{\ast}$, where $(e_{i})$ is any basis of
$V$. Then $K_{l}(A)$ is the $N$-complex
\begin{equation}
\cdots \rightarrow A\otimes A_{n}^{!\ast} \rightarrow A\otimes A_{n-1}^{!\ast}
\rightarrow\cdots \rightarrow A\otimes A_{2}^{!\ast} \rightarrow
A\otimes V \rightarrow A\,,
\end{equation}
where the differential $d:A\otimes A_{n}^{!\ast} \rightarrow A\otimes A_{n-1}^{!\ast}$ is defined
by
$$d(a\otimes \alpha)= (a\otimes \alpha).\xi_{l},\ a\in A,\ \alpha \in A_{n}^{!\ast}. $$
The action of $\xi_{l}$ on $A\otimes A_{n}^{!\ast}$ is defined by
$$(a\otimes \alpha).\xi_{l}= \sum_{i} a e_{i} \otimes \alpha .e_{i}^{\ast},$$ 
and $\alpha . e_{i}^{\ast}$ is the element of $A_{n-1}^{!\ast}$ defined by $\langle
\alpha .e_{i}^{\ast},f\rangle =\langle \alpha , e_{i}^{\ast}f\rangle$, $f$ in $A_{n-1}^{!}$. Here
$ae_{i}$ is the product in $A$, $e_{i}^{\ast}f$ is the product in $A^{!}$, and
$\langle\ ,\,\rangle$ is the natural pairing between a vector space and its dual. More generally,
the algebra
$A\otimes A^{!}$ acts on $A\otimes A^{!\ast}$, where
$A^{!\ast}$ is the graded dual. We have $d^{N}=0$ because $\xi_{l}^{N}=0$ in the algebra $A\otimes
A^{!}$.

There are several manners to contract (3.1) into 2-complexes, by putting together alternately $p$
and $N-p$ arrows $d$, $1\leq p \leq N-1$~\cite{bdvw:homog}. The adequate contraction,
denoted by $K'_{l}(A)$, is obtained from (3.1) by keeping the arrow at the far right,
putting together the $N-1$ previous ones, and continuing alternately :
\begin{equation}
\cdots \rightarrow A\otimes A_{2N}^{!\ast} \rightarrow A\otimes A_{N+1}^{!\ast}
\rightarrow A\otimes A_{N}^{!\ast} \rightarrow
A\otimes V \rightarrow A\,.
\end{equation}
It is easily proved that $K'_{l}(A)$ is isomorphic to the Koszul 2-complex $K$
introduced in Section 2. 

From now on, $A$ is supposed to be generalized Koszul. Our aim is to
make explicit the Yoneda algebra
(or Ext-algebra) $E(A)=$ Ext$_{A}^{\ast}(k,k)$. This algebra is graded by the $E(A)_{i}=$
Ext$_{A}^{i}(k,k)$, $i\geq 0$. In general, its product (the Yoneda product) $\rho :
E(A)\times E(A) \rightarrow E(A)$ can be defined as follows (2.7 in~\cite{ben:rep}, paragraph 7
in~\cite{bou:hom}). Let $(P,\delta)\stackrel{\epsilon}{\longrightarrow}k$ be a projective resolution of $k$
in $A$-grMod. An element of $E(A)$ is the class $[f]$ of a cycle $f$ of the complex
Hom$_{A}(P,k)$. Then $\rho ([f],[g])=[f\circ \tilde{g}]$, where $\tilde{g}$ is any cycle of the
complex Hom$_{A}(P,P)$ such that $\epsilon \circ \tilde{g} =g$. Recall also that Hom$_{A}(P,P)$
is $\mathbf{Z}$-graded by
$$\mathrm{Hom}_{A}(P,P)_{n}= \bigoplus_{i}\mathrm{Hom}_{A}(P_{i},P_{i+n}),$$
with differential $\tilde{\delta}_{n}:\mathrm{Hom}_{A}(P,P)_{n} \rightarrow
\mathrm{Hom}_{A}(P,P)_{n-1}$ defined by the graded commutator $\tilde{\delta}\varphi=
[\delta,\varphi]$. In other words, one has
$$\tilde{\delta}_{n}\varphi(x)=\delta_{i+n}(\varphi(x))-(-1)^{n}\varphi(\delta_{i}(x)),\ \varphi \in
\mathrm{Hom}_{A}(P_{i},P_{i+n}),\ x\in P_{i}.$$
In our situation, take for $P\stackrel{\epsilon}{\longrightarrow}k$ the Koszul resolution
$K'_{l}(A)$ with the usual $\epsilon :A\rightarrow k$. So $P_{i}=A\otimes
A_{n(i)}^{!\ast}$, $i\geq 0$. Since the graded $A$-module $P_{i}$ is generated in degree $n(i)$,
the vector space Hom$_{A}(P_{i},k)$ is graded and concentrated in degree $-n(i)$ ($k$ is
concentrated in degree 0). Therefore Hom$_{A}(P,k)$ is a complex of $k$-grMod whose differential
vanishes. Using the canonical isomorphisms
$$\mathrm{Hom}_{A}(P_{i},k) \cong \mathrm{Hom}_{k}(A_{n(i)}^{!\ast},k) \cong A_{n(i)}^{!},$$
we identify Hom$_{A}(P_{i},k)$ to $A_{n(i)}^{!}$ concentrated in degree $-n(i)$. Thus we get
$$E(A)_{i}\cong A_{n(i)}^{!}, \ i\geq 0.$$

It remains to determine the Yoneda product in this identification. Fix $f\in A_{n(i)}^{!}$
and $g\in A_{n(j)}^{!}$. Denote their Yoneda product by $f\bullet g \in A_{n(i+j)}^{!}$ which has
not to be confused with their product $fg$ in the algebra $A^{!}$. In order to apply the above
definition of the Yoneda product, come back to $f$ in Hom$_{A}(P_{i},k)$ and
$g$ in Hom$_{A}(P_{j},k)$. We have
$$(\mathrm{im}(P_{j+1}\stackrel{\delta_{j+1}}{\longrightarrow} P_{j}))_{n(j+1)} \subseteq (A_{1}\,
\mathrm{or}\, A_{N-1})\otimes A_{n(j)}^{!\ast},$$
according to $n(j+1)$ is $n(j)+1$ or $n(j)+N-1$. For any $v$ in $A_{n(j)}^{!\ast}$, $g(a\otimes
v)=\epsilon(a) g(v)$ vanishes if $a$ is in $A_{1}$ or $A_{N-1}$. Therefore $g\circ \delta_{j+1}=0$.
Denote by $Q'$ the following projective complex of $A$-grMod :
$$\cdots \stackrel{(-1)^{j}\delta}{\longrightarrow} P_{j+2}(n(j))
\stackrel{(-1)^{j}\delta}{\longrightarrow} P_{j+1}(n(j))
\stackrel{(-1)^{j}\delta}{\longrightarrow} P_{j}(n(j)) \stackrel{(-1)^{j}\delta}{\longrightarrow} \cdots ,$$ 
in which the last written term is the 0-degree term of the complex. The shift by
$n(j)$ makes $g : P_{j}(n(j)) \rightarrow k$ of degree 0, so that we can consider the morphism of
complexes $Q'\stackrel{g}{\longrightarrow}k$ in
$A$-grMod where $k$ is considered as a complex concentrated in degree 0. By comparison of this
morphism with the Koszul resolution
$P\stackrel{\epsilon}{\longrightarrow}k$, there exists a morphism of complexes $\tilde{g}:
Q'\rightarrow P$ of $A$-grMod such that $\epsilon \circ \tilde{g}=g$. Clearly $\tilde{g}$ may be
viewed in Hom$_{A}(P,P)$, of degree $-j$. Note that 
$\tilde{g}_{-1}=\tilde{g}_{-2}=\cdots =0$. On the other hand, the sign $(-1)^{j}$ in $Q'$ implies
that
$\tilde{g}$ is a \emph{cycle} of Hom$_{A}(P,P)$. So $f\bullet g$ is identified to the element
$f\circ \tilde{g}_{i}$ of Hom$_{A}(P_{i+j},k)$.

A first consequence is the following. The map $\tilde{g}_{i}: P_{i+j}(n(j))\rightarrow P_{i}$ is
of degree 0. Since $P_{i+j}(n(j))$ is generated in degree $n(i+j)-n(j)$ and $P_{i}$ in degree
$n(i)$, we see that $f\circ \tilde{g}_{i}=0$ if $n(i+j)\neq
n(i)+n(j)$, i.e., if $N>2$ with $i$ and $j$ odd, so that $f\bullet
g=0$ in this case. That is really different from $fg$ !

Assume now $n(i+j) = n(i)+n(j)$, i.e., $N=2$ or ($N>2$ with $i$ or $j$ even). We are going to
give an explicit $\tilde{g}$ by contraction of a natural $N$-complex morphism. Let $Q$ be the
following $N$-complex whose $Q'$ is a contraction (in order to
simplify the notations, the internal shift by $n(j)$ in each term is omitted) :
$$\stackrel{(-1)^{j}d}{\longrightarrow}  A\otimes A_{n(j)+N}^{!\ast} \stackrel{(-1)^{j}d}{\longrightarrow} A\otimes
A_{n(j)+N-1}^{!\ast} \stackrel{d}{\rightarrow} \cdots 
\stackrel{d}{\rightarrow} A\otimes A_{n(j)+1}^{!\ast}
\stackrel{(-1)^{j}d}{\longrightarrow} A\otimes A_{n(j)}^{!\ast}
\stackrel{(-1)^{j}d}{\longrightarrow} $$ 
in which the last written term is the 0-degree term of the
$N$-complex. More explicitly, each term $A_{n(j)+m}^{!\ast}$ for $m \equiv 0\, (\mathrm{mod} N)$ is
surrounded with arrows $(-1)^{j}d$, and the other arrows are $d$.

Let us define an $N$-complex morphism $\tilde{G} :Q \rightarrow K_{l}(A)$ in $A$-grMod as follows. For each $m\geq
0$, $\tilde{G}_{m}=1_{A}\otimes G_{m}$ where $G_{m}:A_{n(j)+m}^{!\ast}\rightarrow A_{m}^{!\ast}$
is defined, for any $\alpha \in A_{n(j)+m}^{!\ast}$, by
$$
\begin{array}{lll}
G_{m} (\alpha) & = g.\alpha & \mathrm{if}\ m \equiv 0\ (\mathrm{mod} N)  \\ 
   & = (-1)^{j} g.\alpha & \mathrm{otherwise}.
\end{array}
$$ 
Here $g$ is seen in $A_{n(j)}^{!}$. The action of $g$ on  
$A_{n(j)+m}^{!\ast}$ is defined by
$$\langle g.\alpha , h \rangle = \langle \alpha , hg \rangle,\ h \in A_{m}^{!},$$ 
where $hg$ is the product in $A^{!}$. The commutativity of the diagram
\begin{eqnarray}
A_{n(j)+m+1}^{!\ast}\ \stackrel{\bar{d}}{\longrightarrow}\ &\ V\otimes A_{n(j)+m}^{!\ast}
\nonumber \\
\downarrow G_{m+1} \ \ \ \ \ \ \     &\ \ \ \ \ \  \downarrow 1_{V}\otimes G_{m}  \\
A_{m+1}^{!\ast}\ \ \ \ \ \ \stackrel{d}{\longrightarrow} & V\otimes A_{m}^{!\ast} \nonumber 
\end{eqnarray}
in which $\bar{d}$ is $d$ or $(-1)^{j}d$, is easy to check. In fact, for any $\alpha$ in $A_{n(j)+m+1}^{!\ast}$,
we have the equalities (with the \emph{same} sign $\pm$) 
$$d\circ G_{m+1}(\alpha)=\pm \sum_{i} e_{i}\otimes ((g.\alpha).e_{i}^{\ast}),$$
$$(1_{V}\otimes G_{m})\circ \bar{d}\,(\alpha)=\pm \sum_{i} e_{i}\otimes (g.(\alpha
.e_{i}^{\ast})),$$
and it is immediate that $(g.\alpha).e_{i}^{\ast}=g.(\alpha.e_{i}^{\ast})$.
 
As $G_{0}=g$, we have $\epsilon \circ \tilde{G}=g$. Next,
$\tilde{g}_{\ell} : P_{\ell+j} \rightarrow P_\ell$ is defined for any $\ell\geq 0$ by
$$
\begin{array}{lll}
\tilde{g}_{\ell} & = \tilde{G}_{n(\ell)} & \mathrm{if}\ \ell \ \mathrm{or}\
j\ \mathrm{even} \\ 
   & =  \tilde{G}_{n(\ell)}\circ d^{N-2} & \mathrm{if}\ \ell \ \mathrm{and}\ j\ \mathrm{odd}.
\end{array} 
$$
Note that $\tilde{g}_{\ell}$ is well-defined in both cases. In the first
case, it is clear because $n(\ell +j)= n(\ell)+n(j)$. In the second case,
since $n(\ell+j)=n(\ell)+n(j)+N-2$, $d^{N-2}: P_{\ell+j} \rightarrow A\otimes
A_{n(\ell)+n(j)}^{!\ast}$ is followed by $\tilde{G}_{n(\ell)} : A\otimes
A_{n(\ell)+n(j)}^{!\ast} \rightarrow P_\ell$. 

Clearly $\tilde{g}:
Q'\rightarrow P$ is a morphism of complexes of $A$-grMod such that $\epsilon \circ \tilde{g}=g$.
So $f\bullet g$ viewed in $A_{n(i+j)}^{!}$ coincides with $f\circ
G_{n(i)}$. Since $n(i) \equiv 0\ (\mathrm{mod} N)$ if and only if $i$ is
even, one has for any $\alpha$
in $A_{n(i+j)}^{!\ast}$,
$$f\circ G_{n(i)}(\alpha)=f((-1)^{ij} g.\alpha)=(-1)^{ij} \langle fg, \alpha \rangle,$$
hence $f\bullet g =(-1)^{ij}fg$. We have obtained :
\Bpo
Let $f$ be in $E(A)_{i}=A_{n(i)}^{!}$ and $g$ in $E(A)_{j}=A_{n(j)}^{!}$. Denote by $f\bullet g$ the
Yoneda product, and $fg$ the product in $A^{!}$. Then we have
$$
\begin{array}{lll}
f\bullet g & = (-1)^{ij}fg & \mathrm{if}\ N=2\ \mathrm{or}\ (N>2\ \mathrm{with}\ i\ \mathrm{or}\
j\ \mathrm{even}) \\ 
   & = 0 & \mathrm{if}\ N>2\ \mathrm{with}\ i\ \mathrm{and}\ j\ \mathrm{odd}.
\end{array} 
$$
\Epo

If $N=2$, we find $f\bullet g =(-1)^{ij}fg$. Note that if $N>2$, the factor $(-1)^{ij}$ is
always $+1$. Show how this factor agrees in any case with the Koszul-Quillen sign rule. The
canonical injection $\mathrm{can} : W_{m} \rightarrow V^{\otimes m}$ has a surjective transposed map
$\mathrm{can}^{\ast} : (V^{\otimes m})^{\ast} \rightarrow W_{m}^{\ast}$ whose kernel is
$I(R^{\perp})_{m}$, once $(V^{\otimes m})^{\ast}$ is identified to $(V^{\ast})^{\otimes m}$.
Write down this identification for $m=2$ (the general case will then be clear) : if
$\alpha$, $\beta$ are in $V^{\ast}$, the element $\alpha
\otimes \beta$ of $V^{\ast}\otimes V^{\ast}$ is identified to the element $a\otimes b \mapsto
\alpha (a) \beta (b)$ of $(V\otimes V)^{\ast}$. It is worth noticing that Beilinson, Ginzburg,
Soergel~\cite{bgso:kdp} adopt another convention by taking $\beta (a) \alpha (b)$ instead of
$\alpha (a) \beta (b)$. But to recover the Koszul-Quillen tensor product of the linear forms as
below with this convention, one should set $E(A)=$ Ext$_{A^{\circ}}^{\ast}(k,k)$ where $A^{\circ}$
is the opposite algebra of $A$. 

From can$^{\ast}$, we deduce the isomorphism $A_{m}^{!} \cong W_{m}^{\ast}$, which sends the class
$\bar{u}$ of an element $u$ of $(V^{\otimes m})^{\ast}$ to its restriction $u|_{W_{m}}$ to
$W_{m}$. Then the product in the algebra $A^{!}$ of $\bar{u} \in A_{m}^{!}$ and $\bar{v} \in
A_{n}^{!}$ is sent to $(u\otimes v)|_{W_{m+n}}$, with our convention on $u \otimes v$
considered in $(V^{\otimes (m+n)})^{\ast}$. If $m=n(i)$ and $n=n(j)$, the Yoneda product
of $\bar{u} \in E(A)_{i}$ and $\bar{v} \in E(A)_{j}$ is sent to $(-1)^{ij} (u\otimes
v)|_{W_{m+n}}$ when $n(i+j)=n(i)+n(j)$. But the Koszul-Quillen tensor product of the linear forms
$u$ and $v$ is also given by
$$\langle u\stackrel{K-Q}{\otimes}v, a\otimes b\rangle = (-1)^{ij} \langle u,a \rangle \langle
v,b \rangle,\ a\in W_{n(i)},\ b\in W_{n(j)},$$
since $v$ has degree $j$, and $a$ has degree $i$.

\setcounter{equation}{0}

\section{Bimodule Koszul $N$-complexes}
Throughout this section, $A$ is an $N$-homogeneous algebra on $V$, with $R$ as space of relations.
We need the two canonical elements $\xi_{l}= \sum_{i} e_{i}\otimes e_{i}^{\ast}$ and $\xi_{r}=
\sum_{i} e_{i}^{\ast}\otimes e_{i}$. In the category $A$-grMod-$A$ of graded
$A$-$A$-bimodules (with morphisms of degree 0), define the following $N$-differentials :
$$d_{l},\ d_{r} :\ \ A\otimes A_{n}^{!\ast}\otimes A \rightarrow A\otimes A_{n-1}^{!\ast}\otimes
A,\ n\geq 0,$$
$$d_{l}(a\otimes \alpha \otimes b)= ((a\otimes \alpha).\xi_{l})\otimes b= \sum_{i} a e_{i} \otimes
\alpha .e_{i}^{\ast} \otimes b,$$
$$d_{r}(a\otimes \alpha \otimes b)= a\otimes (\xi_{r}.(\alpha \otimes b)) = \sum_{i} a\otimes
e_{i}^{\ast}.\alpha \otimes e_{i}b.$$
The left and right actions $\alpha.e_{i}^{\ast}$, $e_{i}^{\ast}.\alpha$ have already been
defined. Since these actions commute, $d_{l}$ and $d_{r}$ commute. Fix a primitive $N$-root of unity
$q$ (we enlarge the ground field $k$ if necessary). Define 
$d:\ A\otimes A_{n}^{!\ast}\otimes A \rightarrow A\otimes A_{n-1}^{!\ast}\otimes A$
by $d=d_{l}-q^{n-1}d_{r}$. Explicitly we have
\begin{equation}
\cdots  \stackrel{d_{l}-d_{r}}\longrightarrow A\otimes A_{N}^{!\ast} \otimes A
\ \stackrel{d_{l}-q^{N-1}d_{r}}\longrightarrow \ \cdots  \stackrel{d_{l}-qd_{r}}\longrightarrow
A\otimes V \otimes A \stackrel{d_{l}-d_{r}}\longrightarrow A\otimes A.
\end{equation}  
Since $\prod_{n=0}^{n=N-1} (d_{l}-q^{n}d_{r})=d_{l}^{N}-d_{r}^{N}=0$, (4.1) is an $N$-complex of
$A$-grMod-$A$, called the \emph{bimodule Koszul $N$-complex} of $A$ and denoted by $K_{l-r}(A)$. 

The \emph{bimodule Koszul 2-complex} $K'_{l-r}(A)$ (see the reference~\cite{rb:nonquad2}
which puts right the faulty definition in~\cite{rb:nonquad}) is the adequate contraction of
$K_{l-r}(A)$,  defined similarly from (4.1) by keeping the arrow at the far right, putting together
the $N-1$ previous ones, and continuing alternately :
\begin{equation}
\cdots  \stackrel{d^{N-1}}\longrightarrow A\otimes A_{N+1}^{!\ast} \otimes A
\stackrel{d}\longrightarrow A\otimes A_{N}^{!\ast} \otimes A \stackrel{d^{N-1}}\longrightarrow
A\otimes V \otimes A \stackrel{d}\longrightarrow A\otimes A.
\end{equation} 
Here $d=d_{l}-d_{r}$ and $d^{N-1}=d_{l}^{N-1}+d_{l}^{N-2}d_{r} + \cdots
+ d_{l}d_{r}^{N-2}+d_{r}^{N-1}$, so that $K'_{l-r}(A)$ makes sense on any ground field. 

Obviously, $K'_{l-r}(A)\otimes_{A} k\cong K'_{l}(A)$ (as 2-complexes). In order to transfer
acyclicity from $K'_{l}(A)$ to $K'_{l-r}(A)$, we need the following result (which also corrects a
faulty version in~\cite{rb:nonquad}). As we shall see, the proof of this result lies upon
well-known ``filtration-gradation" techniques~\cite{nvo:grad}.
\Bpo
Let $A$ be a connected graded algebra. Assume that the complex of $A$-grMod
\begin{equation} \label{comp}
L \stackrel{f}\longrightarrow M \stackrel{g}\longrightarrow N
\end{equation}
is formed of graded-free modules, with $L$ bounded below. Then this complex is exact if the
following is exact :
\begin{equation} \label{rcomp}
k\otimes_{A} L \stackrel{1_{k}\otimes_{A} f}\longrightarrow k\otimes_{A} M
\stackrel{1_{k}\otimes_{A} g}\longrightarrow k\otimes_{A} N.
\end{equation} 
\Epo
\Bdm
Write down the graded-free modules of (\ref{comp}) as $L=A\otimes E$, $M=A\otimes F$,
$N=A\otimes G$, where $E$, $F$, $G$ are graded $k$-vector spaces (tensor products over $k$ are
denoted without subscript). As $L$ is bounded below, there exists an integer $s$ such that
$E_{i}=0$ if $i<s$. Consider $A$ as a filtered $k$-algebra by the 
$$A_{\geq p}=A_{p}\oplus A_{p+1} \oplus \cdots,\ p\geq 0.$$
The filtration is non-increasing, exhaustive, separated. Then $M$ is a filtered $A$-module by the
$A_{\geq p}\otimes F$, $p\geq 0$, and this filtration is non-increasing, exhaustive, separated. Idem for
$L$ and $N$. Since $A_{\geq p}= \bigoplus_{n} (A_{\geq p}\cap A_{n})$, the filtration of $A$ is
compatible with its gradation. It is straightforward to verify that
$$\bigoplus_{n\in \mathbf{Z}}\, (A_{\geq p}\otimes F)\cap (A\otimes F)_{n} = A_{\geq p}\otimes F,$$
which means that the filtration of $M$ is compatible with its gradation. We have the same with $L$
and $N$. Note that $(A_{\geq p}\otimes E)\cap (A\otimes E)_{n}=0$ if $p>n-s$. Therefore the
filtration of $L$ induces on each homogeneous component $L_{n}$ a \emph{finite} (hence complete)
filtration. We can say that the filtration of $L$ is graded-complete.

The complex (\ref{rcomp}) is in $k$-grMod. Identifying naturally $k \otimes_{A} L$ to $E$, ...,
one sees that (\ref{rcomp}) is identified to the $k$-grMod complex
\begin{equation} \label{icomp}
E \stackrel{\tilde{f}}\longrightarrow F \stackrel{\tilde{g}}\longrightarrow G.
\end{equation}
Here $\tilde{f}$ is defined as the composite 
$E\hookrightarrow L \stackrel{f}\longrightarrow M \twoheadrightarrow F,$
where $E\hookrightarrow L$ (resp. $M \twoheadrightarrow F$) is the
canonical injection (resp. projection). Idem for $\tilde{g}$.

Since $A$ is $\mathbf{N}$-graded, $f$ sends $E$ to $A_{\geq 0}\otimes F$. Thus by linearity, $f$
sends $A_{\geq p}\otimes E$ to $A_{\geq p}\otimes F$. So $f$ is compatible with the filtrations of
$L$ and $M$. Similarly $g$ is compatible with the filtrations of $M$ and $N$. By $A_{\geq p}/
A_{\geq p+1} \cong A_{p}$, the associated graded algebra gr$(A)$ is naturally isomorphic to $A$,
and gr$(L)\cong A\otimes E$, gr$(M)\cong A\otimes F$, gr$(N)\cong A\otimes G$. Moreover, gr$(f)$ 
is identified to $1_{A}\otimes \tilde{f} : A\otimes E \rightarrow A\otimes F$. Similarly, gr$(g)$ 
is identified to $1_{A}\otimes \tilde{g} : A\otimes F \rightarrow A\otimes G$.

Assume (\ref{rcomp}) exact. Applying the exact functor $A\otimes -$ to the exact complex
(\ref{icomp}) and using the above identifications, one sees that the complex 
$$\mathrm{gr}(L) \stackrel{\mathrm{gr}(f)}\longrightarrow \mathrm{gr}(M)
\stackrel{\mathrm{gr}(g)}\longrightarrow \mathrm{gr}(N)$$ 
is exact. One concludes thanks to the following proposition, which is the graded analogue of a
classical result (Theorem III.3 in~\cite{nvo:grad}) and whose proof is left to the reader.
\Edm
\Bpo
Let $R$ be a $k$-algebra endowed with a non-increasing $\mathbf{Z}$-filtration $F$. Let $L$, $M$,
$N$ be three filtered $R$-modules. The filtration of these modules is also denoted by $F$. Assume
that the $k$-algebra $R$ is $\mathbf{Z}$-graded and that the $R$-modules $L$, $M$, $N$ are
$\mathbf{Z}$-graded. Assume that the filtrations and the gradations are compatible. Consider the
following complex of filtered $R$-modules
\begin{equation} \label{fcomp}
L \stackrel{f}\longrightarrow M \stackrel{g}\longrightarrow N,
\end{equation}
such that $f$ and $g$ are homogeneous of degree 0 (for the gradations). If the filtration of $L$
is complete in $R$-grMod (i.e., graded-complete) and if the filtration of $M$ is exhaustive and
separated, the exactness of the gr$_{F}(R)$-grMod complex
$$\mathrm{gr}_{F}(L) \stackrel{\mathrm{gr}_{F}(f)}\longrightarrow \mathrm{gr}_{F}(M)
\stackrel{\mathrm{gr}_{F}(g)}\longrightarrow \mathrm{gr}_{F}(N)$$
implies the exactness of (\ref{fcomp}). 
\Epo

Denote by $\mu : A\otimes A \rightarrow A$ the product of $A$. For $\alpha \in V$, we have
$$d_{l}(a\otimes \alpha \otimes b)= \sum_{i} \langle \alpha, e_{i}^{\ast}\rangle a e_{i} \otimes
b,$$
$$d_{r}(a\otimes \alpha \otimes b)= \sum_{i} \langle \alpha, e_{i}^{\ast}\rangle a \otimes e_{i}
b,$$ hence $\mu \circ (d_{l}-d_{r})=0$. Since the beginning $A\otimes R \rightarrow A\otimes V
\rightarrow A \stackrel{\epsilon}\longrightarrow k \rightarrow 0$ of the Koszul 2-complex (with
the augmentation $\epsilon$) is always exact, Proposition 4.1 (or rather, its right modules
version) implies immediately :
\Bpo
The following complex is exact
$$A\otimes A_{N}^{!\ast} \otimes A \stackrel{d^{N-1}}\longrightarrow
A\otimes V \otimes A \stackrel{d}\longrightarrow A\otimes A \stackrel{\mu}\longrightarrow A
\rightarrow 0.$$
\Epo

We can now state the bimodule characterization for generalized Koszulity.
\Bte
The $N$-homogeneous algebra $A$ is Koszul if and only if the following complex is exact :
\begin{equation} \label{bcomp}
K'_{l-r}(A)\stackrel{\mu}\longrightarrow A \rightarrow 0.
\end{equation}
\Ete
\Bdm
Assume $A$ Koszul. Applying the functor $-\otimes_{A} k$ to (\ref{bcomp}), we get the Koszul
2-resolution of the trivial module $k$ in $A$-grMod, and we conclude by Proposition 4.1. Conversely,
assume (\ref{bcomp}) exact. Viewing (\ref{bcomp}) as a projective resolution of $A$ in grMod-$A$ and
comparing it with the identity map of $A$, one sees that there are in grMod-$A$ two morphisms of
resolutions $f: K'_{l-r}(A)
\rightarrow A$ and $ g : A \rightarrow K'_{l-r}(A)$ (where $A$ is considered as a complex
concentrated in degree 0), and a chain homotopy $s: K'_{l-r}(A)
\rightarrow K'_{l-r}(A)$ such
that $$1_{K'_{l-r}(A)}-g f =sd + ds.$$ 
Clearly, $1_{K'_{l-r}(A)}=sd + ds$ in degree $>0$. Applying
$-\otimes_{A} k$, we draw 
$$1_{K'_{l}(A)}=(s\otimes_{A} 1_{k})d_{l} + d_{l}(s\otimes_{A} 1_{k})$$
in degree $>0$, which implies that $K'_{l}(A)$ is exact in degree $>0$.
\Edm

Assume $A$ Koszul. Each bimodule of the projective resolution (\ref{bcomp}) is generated in only
one degree, and these degrees are increasing. Thus the resolution is minimal in
$A$-grMod-$A$. It is called the \emph{bimodule Koszul resolution} of $A$. By
minimality, its length (in $\mathbf{N}\cup \{ \infty \}$) is the Hochschild dimension of
$A$ (the Hochschild dimension of
$A$ is defined as being the projective dimension of $A$ in $A$-grMod-$A$, see just below).
Joining this with Corollary 2.6, we get : 
\Bte
Assume that the $N$-homogeneous algebra $A$ is Koszul. The Hochschild dimension of $A$ coincides
with the global dimension of $A$.
\Ete

For the convenience of the reader, let us recall some basic facts on Hochschild dimension (simply
called dimension by Cartan and Eilenberg~\cite{ce:ha}). Let $k$ be a commutative ring and let $A$
be an associative $k$-algebra with unit. Set $A^{e}=A\otimes A^{\circ}$ (tensor product over $k$). Each
$A$-$A$-bimodule $M$ is a left (right) $A^{e}$-module $M_{l}$ ($M_{r}$) for the action $(a\otimes
b).m = amb$ ($m.(a\otimes b)=bma$). Since $(A^{e})^{\circ}$ is isomorphic to $A^{e}$ by the flip,
it is clear that
$\mathrm{pd}_{A^{e}}M_{l}=\mathrm{pd}_{A^{e}}M_{r}$, which will be denoted by $\mathrm{pd}_{A^{e}}M$.
\Bdf
$\mathrm{pd}_{A^{e}}A$ is called the Hochschild dimension of $A$.
\Edf

Using $Ext_{A^{e}}^{i}(A,M)=HH^{i}(A,M)$, the
terminology comes from 
$$\mathrm{pd}_{A^{e}}A = \sup \{i\in \mathbf{N}\cup \{\infty\}\,;\mathrm{there\ exists}\ M \
\mathrm{such\ that}\ HH^{i}(A,M)\neq 0\,\},$$
with the convention $HH^{\infty}(A,M)=0$. If $k$ is a field, one has inequalities
$$\mathrm{l.gl.}\dim A \leq \mathrm{pd}_{A^{e}}A,\ \ \mathrm{r.gl.}\dim A \leq \mathrm{pd}_{A^{e}}A,$$
which turn out to be equalities in some important circumstances, as group algebras, or enveloping
algebras of Lie algebras~\cite{ce:ha}. So Theorem 4.5 provides another circumstance for equalities.
On the other hand, suppose that $k$ is a field and $A$ is connected graded. For any bounded
below left $A$-module $M$, $\mathrm{pd}_{A}M$ coincides with $\mathrm{gr.pd}_{A}M$ (i.e., the projective dimension of
$M$ computed in the category $A$-grMod), and thus coincides with
the length of a minimal projective resolution of $M$ in $A$-grMod~\cite{nvo:grad}. In particular,
this holds for the algebra $A^{e}$ instead of $A$, and the bimodule $A$ instead of $M$. We have
also used the well-known fact that the (left and right) global
dimension of $A$ is equal to $\mathrm{pd}_{A}(_{A}k)$~\cite{nvo:grad}.

\setcounter{equation}{0}

\section{AS-Gorenstein algebras}
\Bdf
Let $A$ be a connected graded $k$-algebra of finite global dimension $D$. One says that $A$ is
AS-Gorenstein if \emph{Ext}$_{A}^{i}(k,A)=0$ for $i \neq D$, and \emph{Ext}$_{A}^{D}(k,A)=k$. 
\Edf

In this text, Ext will be always understood as a functor of the category $A$-Mod of left
$A$-modules. So AS-Gorenstein means left AS-Gorenstein in this definition. Actually, dealing with
finite global dimension algebras, right AS-Gorenstein is equivalent to left AS-Gorenstein
(\cite{sz:gro}, Proposition 3.1). We have a first link with Koszul property in low dimension. 
\Bpo
Let $A$ be a connected graded $k$-algebra. The space $V$ of generators is concentrated in
degree $1$, and the space $R$ of relations lives in degrees $\geq 2$. Assume that the global dimension
$D$ of $A$ is $2$ or $3$, and that $A$ is AS-Gorenstein. Then $A$ is $N$-homogeneous and Koszul,
with $N=2$ if $D=2$, and $N\geq 2$ if $D=3$.
\Epo
\Bdm
We know that Ext$_{A}^{0}(k,k)=k$ is concentrated in degree 0, Ext$_{A}^{1}(k;k)=V^{\ast}$ is
concentrated in degree $-1$, Ext$_{A}^{2}(k,k)=R^{\ast}$ lives in degrees $\leq -2$. Moreover
Ext$_{A}^{D}(k,k)=k$ is necessarily concentrated in one degree, which will be denoted by $-e$.
Proposition 2.1 shows that
$e\geq D$. AS-Gorenstein property implies that there exist isomorphisms of graded $k$-vector
spaces~\cite{smi:some} 
\begin{equation} \label{isom}
\mathrm{Ext}_{A}^{i}(k,k) \cong (\mathrm{Ext}_{A}^{D-i}(k(e),k))^{\ast},\ 0\leq i \leq D.
\end{equation} 
Assume $D=2$. Then (\ref{isom}) for $i=1$ leads to $e=2$, so that $R$ is concentrated in degree 2
with dim$R=1$. Koszulity is clear. Assume now $D=3$. Then (\ref{isom}) for $i=1$ implies that $R$
is concentrated in degree $e-1$, with $\dim R= \dim V$. Hence $A$ is $N$-homogeneous, $N=e-1$.
Since Ext$_{A}^{3}(k,k)$ is concentrated in degree $-N-1$, $A$ is Koszul.
\Edm

This result is implicit in the paper of Artin and
Schelter~\cite{as:regular}, and it is a
starting point for their classification. Examples for $D=2$ or 3 with $N=2$, and for $D=N=3$, are
in~\cite{as:regular}. Examples in characteristic 0 for $D=3$ with any $N>2$ are given by Theorem 1.1
when $q=1$. Note that if $D=4$ in the hypotheses of the proposition, $A$ can be
non-homogeneous (\cite{as:regular}, Proposition 1.20), but if $A$ is $N$-homogeneous, then $A$ is Koszul (use again
(\ref{isom})) and $N=2$ by Proposition 5.3 below.

We want to investigate AS-Gorenstein property for generalized Koszul algebras. If
$A$ is a Koszul $N$-homogeneous algebra, then the right $A$-module
Ext$_{A}^{i}(k,A)$ is naturally identified to the right $A$-module
$H^{i}(\mathrm{Hom}_{A}(K'_{l}(A),A))$, $i\geq 0$. As Manin showed it in the quadratic
case~\cite{man:kalg}, there is a nice way to express the grMod-$A$ complex
$\mathrm{Hom}_{A}(K'_{l}(A),A)$~\cite{bdvw:homog}. Conformally to our point of view, we begin to
describe an $N$-complex. For any $N$-homogeneous algebra $A$, we introduce the $N$-complex
$L_{r}(A)$ as being $A^{!}\otimes A$ graded by the $A_{n}^{!}\otimes A$, $n\geq 0$, and endowed with
the $N$-differential $\xi_{r}\cdot$ (left multiplication by $\xi_{r}= \sum_{i} e_{i}^{\ast}\otimes
e_{i}$). Viewing $A_{n}^{!}$ as concentrated in degree $-n$, $L_{r}(A)$ is a cochain $N$-complex of
grMod-$A$. Then, it is easy to find a natural isomorphism of cochain $N$-complexes of grMod-$A$
between $L_{r}(A)$ and $\mathrm{Hom}_{A}(K_{l}(A),A)$~\cite{bdvw:homog}. The $N$-complex
$L_{l}(A)$ is defined analogously.

Let us write down $L_{r}(A)$ more explicitly :
\begin{equation} \label{lr}
A \stackrel{\xi_{r}\cdot} \longrightarrow V^{\ast} \otimes A \stackrel{\xi_{r}\cdot} \longrightarrow
\cdots
\stackrel{\xi_{r}\cdot}\longrightarrow A_{n}^{!}\otimes A \stackrel{\xi_{r}\cdot}\longrightarrow
A_{n+1}^{!}\otimes A \stackrel{\xi_{r}\cdot}\longrightarrow \cdots 
\end{equation}
If we keep the first arrow, put together the $N-1$ following ones, and continue
alternately, we define the adequate contraction $L'_{r}(A)$ :
\begin{equation} \label{clr}
A \stackrel{\xi_{r}\cdot} \longrightarrow V^{\ast} \otimes A \ \stackrel{\xi_{r}^{N-1}\cdot}
\longrightarrow \  A_{N}^{!}\otimes A \stackrel{\xi_{r}\cdot}\longrightarrow
A_{N+1}^{!}\otimes A \ \stackrel{\xi_{r}^{N-1}\cdot}\longrightarrow \ \cdots 
\end{equation}
Then the grMod-$A$ 2-complexes $L'_{r}(A)$ and $\mathrm{Hom}_{A}(K'_{l}(A),A)$ are naturally
isomorphic. Comparing $L'_{r}(A)$ with $K'_{r}(A)$, we are going to give a convenient criterion for
generalized Koszul algebras to be AS-Gorenstein.

We work in an intermediate class of algebras formed of the Koszul $N$-homogeneous algebras $A$
with global dimension $D<\infty$ and Ext$_{A}^{D}(k,A)=k$. Let $A$ be such an algebra. Since
Ext$_{A}^{D}(k,A)= H^{D}(L'_{r}(A))$, there is a linear map $\epsilon' : A_{n(D)}^{!}\otimes A
\rightarrow k$ such that the small complex
$$\longrightarrow A_{n(D)}^{!}\otimes A \stackrel{\epsilon'} \longrightarrow k \longrightarrow 0,$$
in which the first arrow is $\xi_{r}\cdot$ or $\xi_{r}^{N-1}\cdot$, is exact. The image of the first
arrow is contained in $A_{n(D)}^{!}\otimes A_{\geq 1}$ or $A_{n(D)}^{!}\otimes A_{\geq N-1}$, so
$\epsilon'$ is injective on $A_{n(D)}^{!}\otimes A_{0}$ or $A_{n(D)}^{!}\otimes (A_{0}\oplus
\cdots \oplus A_{N-2})$. Therefore, dim$(A_{n(D)}^{!})=1$ and, if $N>2$, the second case does not
occur i.e., $D$ is odd. We have proved :
\Bpo
Assume that $A$ is a Koszul $N$-homogeneous algebra with global dimension $D<\infty$ and
$\mathrm{Ext}_{A}^{D}(k,A)=k$. Then $\mathrm{dim}(A_{n(D)}^{!})=1$ and, if $N>2$, $D$ is odd.
Moreover the following additivity holds
$$n(i)+n(D-i)=n(D), \ 0\leq i \leq D.$$
\Epo

Let us keep $A$ in our intermediate class of algebras. Since $A$ is Koszul of global dimension
$D$, we have $A_{n(D)+N-1}^{!\ast}=0$, hence $A_{i}^{!}=0$ for $i\geq n(D)+N-1$ (Warning!
$A_{i}^{!}\neq 0$ can occur if $n(D)<i<n(D)+N-1$ : take $N>2$, $R=0$ and dim$V=1$). So the
cochain complex $L'_{r}(A)$ is the following
\begin{equation} \label{eclr}
A \stackrel{\xi_{r}\cdot} \longrightarrow V^{\ast} \otimes A \stackrel{\xi_{r}^{N-1}\cdot}
\longrightarrow A_{N}^{!}\otimes A \stackrel{\xi_{r}\cdot}\longrightarrow \cdots 
\stackrel{\xi_{r}^{N-1}\cdot} \longrightarrow 
A_{n(D-1)}^{!}\otimes A \stackrel{\xi_{r}\cdot}\longrightarrow A_{n(D)}^{!}\otimes A \rightarrow 0. 
\end{equation}
Denote by $L'_{r}(A)^{ch}$ the naturally associated chain complex, i.e., in (\ref{eclr}) the
complex degrees are successively $0, -1, -2, \ldots, -D+1, -D$. For comparing with the chain complex
$K'_{r}(A)$, we need the chain complex $L'_{r}(A)^{ch}[-D]$ where the successive complex degrees are
$D,D-1,D-2, \ldots, 1,0$. It is clear that $L'_{r}(A)^{ch}[-D]$ is also the adequate contraction
$L_{r}(A)^{ch}[-n(D)]'$. Here $[\ \ ]$ denotes the shift of $N$-complexes.

Fix a generator $u$ of $A_{n(D)}^{!}$ and choose $\epsilon' = \epsilon_{u}$, where
$\epsilon_{u}$ is defined by $\epsilon_{u}(u\otimes 1)=1$. Denoting as usual by
$(\ \ )$ the shift in grMod-$A$, $\epsilon_{u}$ is a morphism of grMod-$A$ complexes of
$L_{r}(A)^{ch}[-n(D)](-n(D))'$ to $k$ (recall that $A_{i}^{!}$ is concentrated in degree $-i$). 
We can also consider $\epsilon_{u}$ as a morphism of grMod-$A$ $N$-complexes from
$L_{r}(A)^{ch}[-n(D)](-n(D))$ to $k$.

Introduce now the left $A^{!}\otimes A$-linear map $\Phi : A^{!}\otimes A \rightarrow
A^{!\ast}\otimes A$ by $\Phi = \bar{\Phi}\otimes 1_{A}$, $\bar{\Phi} : A^{!} \rightarrow
A^{!\ast}$, and $\bar{\Phi}(1)=u^{\ast}$. Here $u^{\ast}$ is the basis of $A_{n(D)}^{!\ast}$ dual
to the basis $u$ of $A_{n(D)}^{!}$. Clearly, $\Phi$ is graded by the
$\Phi_{i}=\bar{\Phi}_{i}\otimes 1_{A}$, $\bar{\Phi}_{i}: A^{!}_{n(D)-i}\rightarrow
A^{!\ast}_{i}$. The $N$-differentials of $L_{r}(A)$ and $K_{r}(A)$ are both left
multiplications by $\xi_{r}$. Therefore, taking into account the degrees, one sees that 
$$\Phi : L_{r}(A)^{ch}[-n(D)](-n(D)) \rightarrow K_{r}(A)$$ 
is a morphism of grMod-$A$ $N$-complexes. Since $\bar{\Phi}_{0}: A_{n(D)}^{!} \rightarrow k$
coincides with $u^{\ast}$, one has $\epsilon \circ \Phi = \epsilon_{u}$. The adequate contraction
$\Phi'$ of $\Phi$ is the following morphism of grMod-$A$ complexes:
$$\Phi' : L_{r}(A)^{ch}[-n(D)](-n(D))' \rightarrow K'_{r}(A),$$ 
formed by $\Phi_{0}, \Phi_{1}, \Phi_{N}, \Phi_{N+1}, \ldots, \Phi_{n(D)-1}, \Phi_{n(D)}$. One has 
$\epsilon \circ \Phi' = \epsilon_{u}$.

Then $A$ is AS-Gorenstein if and only if $H^{i}(L'_{r}(A))=0$ for any $i>0$, i.e., $\epsilon_{u} :
L_{r}(A)^{ch}[-n(D)](-n(D))' \rightarrow k$ is a resolution, necessarily minimal projective in
grMod-$A$ (minimality comes from the fact that each module of the projective resolution is
generated in only one degree, and the degrees are increasing).
Knowing that a morphism between two minimal projective resolutions is an isomorphism, we have
obtained the following criterion.
\Bte
Assume that $A$ is a Koszul $N$-homogeneous algebra with global dimension $D<\infty$ and
$\mathrm{Ext}_{A}^{D}(k,A)=k$. Fix a generator $u$ of the vector space $A_{n(D)}^{!}$. For
$0\leq i \leq D$, define the linear maps $\bar{\Phi}_{n(i)}: A^{!}_{n(D-i)}\rightarrow
A^{!\ast}_{n(i)}$ by $f\mapsto f.u^{\ast}$. Then $A$ is AS-Gorenstein if and only if these maps
are bijective.
\Ete

\addtocounter{Df}{1}
\Rm
$\bar{\Phi}_{0}: u\mapsto 1$ and $\bar{\Phi}_{n(D)}: 1\mapsto u^{\ast}$ are bijective. On the
other hand, the commutative diagram with exact lines
\begin{eqnarray}
A^{!}_{n(D-1)}\otimes A \stackrel{\xi \cdot}{\longrightarrow} & A^{!}_{n(D)}\otimes A
\stackrel{\epsilon_{u}}{\longrightarrow} & k
\nonumber \\
\downarrow \Phi_{1} \ \ \ \ \   & \downarrow \Phi_{0}\ \ \  & \downarrow 1_{k}  \\
V\otimes A \ \stackrel{d}{\longrightarrow} & A
\ \ \stackrel{\epsilon}{\longrightarrow} & k \nonumber 
\end{eqnarray}
shows that $d(\Phi_{1}(A^{!}_{n(D-1)}\otimes A))=A_{\geq 1}= d(V\otimes A)$. As the map
$d:V\otimes A \rightarrow d(V\otimes A)$ is essential, we have $\Phi_{1}(A^{!}_{n(D-1)}\otimes
A)=V\otimes A$. Thus $\bar{\Phi}_{1}$ is surjective.

\addtocounter{Df}{1}
\Rm
Since $\Phi_{n(D)}$ and $d:A^{!\ast}_{n(D)}\otimes A \rightarrow A^{!\ast}_{n(D-1)}\otimes A$ are
injective (the second one by Koszulity), $\xi \cdot : A\rightarrow A^{!}_{1}\otimes A$ is
injective. Thus Hom$_{A}(k,A)=0$.

\addtocounter{Df}{1}
\Rm
$\bar{\Phi} : A^{!} \rightarrow A^{!\ast}$ is bijective if and only if
the finite dimensional algebra $A^{!}$ is Frobenius and the highest
$i$ such that $A^{!}_{i}\neq 0$ is $i=n(D)$ (see~\cite{smi:some} for a
detailed account on graded Frobenius algebras). If we take $N>2$,
$\dim V =1$, $R=0$, then $A$ is $N$-homogeneous, Koszul, AS-Gorenstein, of
global dimension 1, and $A^{!}$ is
Frobenius with $\bar{\Phi}$ non bijective. Note also that if $A^{!}$
is Frobenius and $i>n(D)$, then $N>2$, $\dim V =1$ and $R=0$.
\Bcr
Assume that $A$ is a Koszul $N$-homogeneous algebra with global dimension $2$.
If $\mathrm{Ext}_{A}^{2}(k,A)=k$, then $A$ is AS-Gorenstein.
\Ecr
\Bdm
Immediate from Theorem 5.4 and Remark 5.5 since $\bar{\Phi}_{0}$, $\bar{\Phi}_{N}$ are bijective,
and $\bar{\Phi}_{1} : V^{\ast} \rightarrow V$ is surjective.
\Edm

\addtocounter{Df}{1}
\Ea
$A$ has $n\geq 2$ generators $x_{1},\ldots ,x_{n}$, and one relation
$$x_{1}x_{2}+x_{2}x_{3}+ \cdots + x_{n-1}x_{n}-x_{n}x_{1}=0.$$
For $n=2$, $A=k[x_{1},x_{2}]$. $A$ is $X$-confluent for the basis $X: x_{1}<x_{2}<\cdots <x_{n}$
since there is no ambiguous word of length 3 (the only nonreduced word of length 2 is
$x_{n}x_{1}$). So $A$ is Koszul. A basis of $A$ is formed of words not containing $x_{n}x_{1}$ as
a factor. A basis of $R^{\perp}$ is formed of words $x_{i}^{\ast}x_{j}^{\ast}$ such that
$x_{i}x_{j}$ is not in the relation of $A$, and of binomials
$x_{i}^{\ast}x_{i+1}^{\ast}+x_{n}^{\ast}x_{1}^{\ast}$, $1\leq i \leq n-1$. The quadratic algebra
$A^{!}$ is $X^{\ast}$-confluent for $X^{\ast} :x_{1}^{\ast}>x_{2}^{\ast}>\cdots >x_{n}^{\ast}$. A
basis of $A_{2}^{!}$ is $x_{n}^{\ast}x_{1}^{\ast}$, and $A_{i}^{!}=0$ for $i>2$. Thus gl.$\dim A
=2$, and the Hilbert series of the graded algebra $A$ is the following:
$$h_{A}(t)=\frac{1}{1-nt+t^{2}}\,.$$
For $n>2$, the latter has a pole between 0 and 1, so GK.$\dim A =\infty$. A
theorem by Stephenson and Zhang~\cite{sz:gro} implies that $A$ is not (left or right)
noetherian in this case. Let us show that $A$ is AS-Gorenstein by using Corollary 5.8. It suffices
to check that the complex 
$$A^{!}_{1}\otimes A \stackrel{\xi \cdot}{\longrightarrow}  A^{!}_{2}\otimes A
\stackrel{\epsilon_{u}}{\longrightarrow}  k$$
is exact. Here $\xi = \sum_{i=1}^{n} x^{\ast}_{i}\otimes x_{i}$ and $u=x_{n}^{\ast}x_{1}^{\ast}$.
From the relations of
$A^{!}$, it is easy to prove that im$(\xi \cdot)$ is the set of the elements 
$$x_{n}^{\ast}x_{1}^{\ast}\otimes (x_{n}a_{1}-x_{1}a_{2}-\cdots -x_{n-1}a_{n}),\ a_{1},\dots ,
a_{n}\ \mathrm{in}\ A.$$
Hence im$(\xi \cdot)= A^{!}_{2}\otimes A_{\geq 1}$ which is ker$(\epsilon_{u})$.

The next corollary is Theorem 1.1 of the introduction. It would be much harder to prove it
by showing directly that $L'_{r}(A)$ is exact. Even if $n=N=3$, the computation of the 
$H^{i}(L'_{r}(A))$ is rather involved.
\Bcr
For $n\geq N\geq 2$, consider the $N$-homogeneous algebra $A$ over $k$ of characteristic 0, with
generators
$x_{1},\ldots ,x_{n}$. The relations of $A$ are the antisymmetrizers of degree $N$, i.e., are the
following
$$\sum_{\sigma \in \mathbf{S}_{N}} \mathrm{sgn} (\sigma)\ x_{i_{\sigma(1)}}\ldots
x_{i_{\sigma(N)}}=0\, ,\ 1\leq i_{1}<\cdots <i_{N} \leq n\,,$$
where $\mathbf{S}_{N}$ denotes the permutation group of $\{1, \ldots ,N\}$. We
know~\cite{rb:nonquad} that $A$ is generalized Koszul, of finite global dimension. Then $A$ is
AS-Gorenstein if and only if either $N=2$ or ($N>2$ \emph{and} $n= Nq + 1$ \emph{for some integer}
$q\geq 1$).
\Ecr
\Bdm
Assume that $A$ is AS-Gorenstein with $N>2$. Using Proposition 5.3, $D$ is odd. 
Write down $D=2q+1$. According to~\cite{rb:nonquad}, $n=Nq+r$ with $1\leq r\leq N-1$. Theorem 5.4 
shows that $A^{!}_{n(D-i)}$ has the same dimension as $A^{!}_{n(i)}$, $0\leq i \leq D$.
Recall~\cite{rb:nonquad} that for $N\leq m \leq n$, $A^{!}_{m}$ has a basis formed of the
words $x_{i_{1}}^{\ast}\cdots x_{i_{m}}^{\ast}$ with $i_{1}>\cdots
>i_{m}$ (and $A^{!}_{m}=0$ if $m>n$). Then
$$ \left(\begin{array}{c}
n \\ n(D)-n(i) 
\end{array} \right)  \ = \ 
\left(\begin{array}{c}
 n \\ n(i) 
\end{array} \right),$$   
which implies $n(D)=2n(i)$ or $n(D)=n$. The first equality does not hold for any $i$, so
$n(D)=n$. As $n(D)=Nq+1$, one has $r=1$. 

Conversely, assume $n=Nq+1$ for some integer $q\geq 1$, and $N\geq 2$ (if $N=2$, the following
proof works for $n$ arbitrary, but $A$ is just a polynomial algebra in this
case!). By~\cite{rb:nonquad},
$D=2q+1$, hence $n(D)=n$. Recall also that a basis of $R^{\perp}$ is formed of words
$x_{j_{1}}^{\ast}\ldots x_{j_{N}}^{\ast}$ having two identical letters, and of binomials
$x_{j_{1}}^{\ast}\ldots x_{j_{N}}^{\ast}- s\, x_{i_{1}}^{\ast}\ldots
x_{i_{N}}^{\ast}$ where $j_{1},\ldots ,j_{N}$ are distinct, not in the decreasing order,
$\{i_{1},\ldots ,i_{N}\}=\{j_{1},\ldots ,j_{N}\}$ with $i_{1}>\cdots >i_{N}$, and $s$ is the sign of
the permutation $i_{1} \mapsto j_{1}, \ldots ,i_{N}\mapsto j_{N}$. First prove
Ext$_{A}^{D}(k,A)=k$, i.e., the complex
$$A^{!}_{n-1}\otimes A \stackrel{\xi \cdot}{\longrightarrow}  A^{!}_{n}\otimes A
\stackrel{\epsilon_{u}}{\longrightarrow}  k$$
is exact. Here $\xi = \sum_{i=1}^{n} x^{\ast}_{i}\otimes x_{i}$ and
$u=x_{n}^{\ast}\cdots x_{1}^{\ast}$. Clearly $\ker(\epsilon_{u})=A^{!}_{n}\otimes A_{\geq 1}$. On
the other hand, for $1\leq i \leq n$, set $u_{i}=x_{n}^{\ast}\cdots \vee_{i} \cdots x_{1}^{\ast}$
where $\vee_{i}$ means that $x_{i}^{\ast}$ is removed. Then the $u_{i}$ form the basis of
$A^{!}_{n-1}$, and we have $(\xi \cdot)(u_{i}\otimes 1) = \pm u\otimes x_{i}$. Therefore
$(\xi \cdot)(A^{!}_{n-1}\otimes A_{0})=A^{!}_{n}\otimes A_{1}$. Since $A_{1}.A= A_{\geq 1}$, we
conclude that $(\xi \cdot)(A^{!}_{n-1}\otimes A)=A^{!}_{n}\otimes A_{\geq 1}$.

Using Theorem 5.4, it remains to prove that $\bar\Phi_{n(i)}$ is bijective for $0\leq i \leq D$.
Actually, we shall prove that $\bar\Phi_{i}$ is bijective for $0\leq i\leq n$ such that $i$ and
$n-i$ are different from $2,3,\ldots ,N-1$. In this case, one has already $\dim A^{!}_{n-i}=\dim
A^{!\ast}_{i}$. For $j_{1}> \cdots >j_{i}$, the words 
$f_{j_{1}\ldots j_{i}}=x_{n}^{\ast}\cdots \vee_{j_{1}} \cdots
\vee_{j_{2}} \cdots \vee_{j_{i}} \cdots x_{1}^{\ast},$
in which the letters $x_{j_{1}}^{\ast},\ldots , x_{j_{i}}^{\ast}$ are removed, form the basis of
$A^{!}_{n-i}$. For $k_{1}>\cdots  >k_{i}$, one has
$$\langle x_{k_{1}}^{\ast}\ldots x_{k_{i}}^{\ast}\,,\, \bar{\Phi}_{i}(f_{j_{1}\ldots j_{i}})\rangle
=\langle x_{k_{1}}^{\ast}\ldots x_{k_{i}}^{\ast}\,f_{j_{1}\ldots j_{i}}\,,\, u^{\ast} \rangle,$$
which amounts to $\pm 1$ if $\{k_{1},\ldots  ,k_{i}\}= \{j_{1},\ldots
,j_{i}\}$, and 0 otherwise. Therefore $\bar{\Phi}_{i}(f_{j_{1}\ldots j_{i}})$ coincides up to a sign with the element
$(x_{j_{1}}^{\ast}\ldots x_{j_{i}}^{\ast})^{\ast}$ of the basis which is dual to the basis of
$A^{!}_{i}$. Thus $\bar{\Phi}_{i}$ is bijective.
\Edm

\addtocounter{Df}{1}
\Rm
If in the previous proof, one has $N>2$, $n=Nq+1$ and $2\leq i \leq N-1$, then $\bar{\Phi}_{i}$ is
not bijective. Hence $A^{!}$ is not Frobenius (Remark 5.7). Note also that for $N>2$, an algebra $A$ as in Corollary 5.10 contains the free
algebra in $x_{1}$ and $x_{2}$, so GK.dim$A=\infty$ and $A$ is not (left or right) noetherian.

The next corollary is Theorem 1.2 of the introduction, generalizing a result by Smith.
Our proof uses the Yoneda product computed in Section 3.
\Bcr
Let $A$ be $N$-homogeneous, generalized Koszul, of finite global dimension $D$. Denote by $E(A)$ the
Yoneda algebra of $A$. Then $A$ is AS-Gorenstein if and only if the
finite dimensional algebra $E(A)$ is Frobenius.
\Ecr
\Bdm
Assume that $A$ is AS-Gorenstein. Choose a
generator $u$ of $A^{!}_{n(D)}$. First we define an arrow $\varphi: E(A) \rightarrow
E(A)^{\ast}(-D)$ in $E(A)$-grMod as follows. For $0\leq i \leq D$, $E(A)_{i}=A^{!}_{n(i)}$. The maps
$\varphi_{i}:E(A)_{i}\rightarrow (E(A)_{D-i})^{\ast}=E(A)^{\ast}(-D)_{i}$ defined by 
$$\varphi_{i}=(-1)^{i.(D-i)} \bar{\Phi}_{n(D-i)}$$
provide an arrow $\varphi: E(A) \rightarrow E(A)^{\ast}(-D)$ in $k$-grMod. It remains to prove the
left $E(A)$-linearity. Here $E(A)^{\ast}$ is naturally a left $E(A)$-module. Denote this action by
$\bullet$ as the Yoneda product. Take $f\in E(A)_{i}$ and $\alpha \in E(A)^{\ast}_{i+j}$. 
Proposition 3.1 shows that $f\bullet \alpha$ vanishes if $N>2$ with $i$ and $j$ odd, and
otherwise, $f\bullet \alpha =(-1)^{i.j} f\cdot \alpha$ where $f\cdot\alpha$ denotes the natural
left action of $A^{!}$ on $A^{!\ast}$. For $f\in E(A)_{i}$ and $g\in E(A)_{j}$, one has to check the
equality $\varphi(f\bullet g)=f\bullet \varphi(g)$. It is clear if $N>2$ with $i$ and $j$ odd, since
$f\bullet g=0= f\bullet \varphi(g)$ in this case (remember that $D$ is odd!). Assume now
$n(i+j)=n(i)+n(j)$. Using Proposition 3.1 and the left $A^{!}$-linearity of $\bar{\Phi}$, we have
$$\varphi(f\bullet g)=(-1)^{ij+(i+j).(D-i-j)} f\cdot \bar{\Phi}(g).$$
We have on the other hand
$$f\bullet \varphi(g)= (-1)^{i.(D-i-j)+j.(D-j)} f\cdot \bar{\Phi}(g).$$
The sign is the same as the previous one, hence the equality (note
that all the signs are $+1$ if $N>2$). Then Theorem 5.4 shows that $\varphi$
is an isomorphism, so that $E(A)$ is Frobenius.

Conversely, assume that the finite dimensional graded algebra $E(A)$ is Frobenius. The highest $i$
such that $E(A)_{i}\neq 0$ is $D$, since $A$ is Koszul of global dimension $D$. According to Lemma
3.2 of~\cite{smi:some}, there is an isomorphism $\psi : E(A) \rightarrow E(A)^{\ast}(-D)$ in
$E(A)$-grMod. It is worth noticing that $D$ is odd if $N>2$. In fact, $\psi(1)$ belongs to
$(E(A)_{1+(D-1)})^{\ast}$, and $\psi(e_{i}^{\ast})=e_{i}^{\ast} \bullet \psi(1)$ vanishes if $D$ is
even. Accordingly $\psi_{i}$ sends $A^{!}_{n(i)}$ to $A^{!\ast}_{n(D)-n(i)}$. So $\psi \otimes
1_{A}$ is a right $A$-linear 0-graded map from the cochain complex $L'_{r}(A)$ to the cochain
complex $K'_{r}(A)^{coch}[-D]$. The next lemma shows that $\psi \otimes
1_{A}$ is a morphism of complexes unless if $N=2$. Checking the lemma
is easier if $N>2$ (no sign occurs!), but we prefer a set-up including
the quadratic case.
\Blm
Setting $\tilde{\psi}_{i}=(-1)^{i.(D-i)} \psi_{i}$ for $0\leq i \leq D$,
$\tilde{\psi}\otimes 1_{A}$ is a morphism of complexes.
\Elm
Assume that the lemma holds. Then $\tilde{\psi}\otimes 1_{A}$ realizes an isomorphism from
the complex $L'_{r}(A)\stackrel{\tilde{\psi}_{D}\otimes \epsilon}\longrightarrow k \rightarrow 0$ to
the complex $K'_{r}(A)^{coch}[-D]\stackrel{\epsilon}\longrightarrow k \rightarrow 0$. The latter
is exact, thus the former is exact. As $L'_{r}(A)$ computes the Ext$_{A}^{i}(k,A)$, $A$ is
AS-Gorenstein. It remains to prove the lemma, i.e., the commutativity of the diagram
\begin{eqnarray}
E(A)_{i}\otimes A\ \ \longrightarrow &\
E(A)_{i+1}\otimes A
\nonumber \\
\downarrow \tilde{\psi}_{i}\otimes 1_{A} \ \ \ \ \ \     &\ \ \ \ \ \ \downarrow 
\tilde{\psi}_{i+1}\otimes 1_{A}
\\ E(A)_{D-i}^{\ast}\otimes A\ \ \longrightarrow&
E(A)_{D-i-1}^{\ast}\otimes A,
\nonumber 
\end{eqnarray}
where the differentials are both $\xi\cdot$ or $\xi^{N-1}\cdot$. Fix $f$ in $E(A)_{i}$. Assume
firstly $N=2$ or ($N>2$ and $i$ even). The differentials are both $\xi\cdot$. Expressing the
action $\cdot$ into the action $\bullet$, we have
$$(\xi\cdot) \circ (\tilde{\psi}_{i}\otimes 1_{A})(f\otimes 1)=(-1)^{i.(D-i)+D-i-1}
\sum_{j}
\psi(e_{j}^{\ast}\bullet f)\otimes e_{j}.$$
On the other hand, we have
$$(\tilde{\psi}_{i+1}\otimes 1_{A})\circ (\xi\cdot)(f\otimes
1)=(-1)^{(i+1).(D-i-1)+i}
\sum_{j}
\psi(e_{j}^{\ast}\bullet f)\otimes e_{j}.$$
Hence the result. Assume now $N>2$ and $i$ odd. The differentials are
both $\xi^{N-1}\cdot$. Since $A^{!}$ is strongly graded, write down $f=\sum_{j} e^{\ast}_{j} f_{j}$
with $f_{j}$ in $E(A)_{i-1}$ (note that $n(i)-1=n(i-1)$). Starting with
$$(\xi^{N-1}\cdot) \circ (\tilde{\psi}_{i}\otimes 1_{A})(f\otimes 1)=
(-1)^{i.(D-i)}\sum_{j_{1}\ldots j_{N-1},j}((e^{\ast}_{j_{1}}\ldots
e^{\ast}_{j_{N-1}})\cdot \psi (e^{\ast}_{j} f_{j})) \otimes e_{j_{1}}\ldots e_{j_{N-1}}$$
the right-hand side is next reduced to
$$(-1)^{i.(D-i)+i-1+D-i}\sum_{j_{1}\ldots j_{N-1}}\psi (e^{\ast}_{j_{1}}\ldots
e^{\ast}_{j_{N-1}} f)\otimes e_{j_{1}}\ldots e_{j_{N-1}}.$$
On the other hand, $(\tilde{\psi}_{i+1}\otimes 1_{A})\circ (\xi^{N-1}\cdot)(f\otimes 1)$ is equal
to 
$$(-1)^{(i+1).(D-i-1)}\sum_{j_{1}\ldots j_{N-1}}\psi
(e^{\ast}_{j_{1}}\ldots e^{\ast}_{j_{N-1}} f)\otimes e_{j_{1}}\ldots e_{j_{N-1}}.$$
Hence the result again.
\Edm

\setcounter{equation}{0}

\section{Poincar\'e duality}
In this section, $A$ is an $N$-homogeneous algebra which is generalized Koszul, AS-Gorenstein and 
of finite global dimension $D$. We want to apply the Van den Bergh duality 
theorem~\cite{vdb:dual} to $A$. Since $A$ has a finite Hochschild dimension (Theorem 4.5), it
remains to prove that Ext$_{A^{e}}^{i}(A, A^{e})=0$ if $i\neq D$, and that
$U=\mathrm{Ext}_{A^{e}}^{D}(A, A^{e})$ is an invertible bimodule. By Theorem 4.4, one has
Ext$_{A^{e}}^{i}(A, A^{e})=H^{i}(\mathrm{Hom}_{A^{e}}(K'_{l-r}(A),A^{e}))$, $i\geq 0$. We are
going to describe the $A$-grMod-$A$ complex $\mathrm{Hom}_{A^{e}}(K'_{l-r}(A),A^{e}))$ by the
bimodule version $L'_{r-l}(A)$ of $L'_{r}(A)$. Here again, we begin to define the $N$-complex
$L_{r-l}(A)$.

For any $N$-homogeneous algebra $A$, set $L_{r-l}(A)=A\otimes A^{!}\otimes A$, graded by the
$A\otimes A_{n}^{!}\otimes A$, $n\geq 0$, and endowed with the cochain $N$-differentials
$\partial_{r}=1_{A}\otimes (\xi_{r}\cdot)$ and $\partial_{l}=(\cdot \xi_{l})\otimes
1_{A}$. These $N$-differentials commute.  Fix a primitive $N$-root of unity
$q$ (and enlarge $k$ if necessary). Define 
$\partial:\ A\otimes A_{n}^{!\ast}\otimes A \rightarrow A\otimes A_{n+1}^{!\ast}\otimes A$ by
$\partial=\partial_{r}-q^{n}\partial_{l}$. Explicitly we have
\begin{equation}
A\otimes A \stackrel{\partial_{r}-\partial_{l}}\longrightarrow A\otimes V^{\ast} \otimes A
\ \stackrel{\partial_{r}-q\partial_{l}}\longrightarrow  \cdots 
\stackrel{\partial_{r}-q^{N-1}\partial_{l}}\longrightarrow A\otimes A_{N}^{!} \otimes A\
\stackrel{\partial_{r}-\partial_{l}}\longrightarrow \cdots.
\end{equation}  
Viewing $A_{n}^{!}$ as concentrated in degree $-n$, $(L_{r-l}(A), \partial)$ is a cochain
$N$-complex of $A$-grMod-$A$. Using the canonical grMod-$A^{e}$ isomorphisms (where $E$ is a finite
dimensional graded $k$-vector space)
$$\mathrm{Hom}_{A^{e}}(A^{e}\otimes E,A^{e}) \cong \mathrm{Hom}_{k}(E,A^{e}) \cong E^{\ast}
\otimes A^{e},$$
and the identifications $A\otimes E \otimes A \cong A^{e}\otimes E$, $E^{\ast}
\otimes A^{e} \cong A\otimes E^{\ast} \otimes A$, one constructs easily an isomorphism of
$A$-grMod-$A$ cochain $N$-complexes from $\mathrm{Hom}_{A^{e}}(K_{l-r}(A),A^{e}))$ to $L_{r-l}(A)$.

If we keep in (6.1) the first arrow, put together the $N-1$ following ones, and continue
alternately, we define the adequate contraction $L'_{r-l}(A)$ :
\begin{equation} 
A\otimes A \stackrel{\partial} \longrightarrow A\otimes V^{\ast} \otimes A \
\stackrel{\partial^{N-1}}
\longrightarrow \  A \otimes A_{N}^{!}\otimes A \stackrel{\partial}\longrightarrow
A\otimes A_{N+1}^{!}\otimes A \ \stackrel{\partial^{N-1}}\longrightarrow \ \cdots 
\end{equation}
Here $\partial=\partial_{r}-\partial_{l}$ and
$\partial^{N-1}=\partial_{r}^{N-1}+\partial_{r}^{N-2}\partial_{l} +
\cdots + \partial_{r}\partial_{l}^{N-2}+\partial_{l}^{N-1}$, so that $L'_{r-l}(A)$ makes sense on
any ground field. Then the $A$-grMod-$A$ 2-complexes
$\mathrm{Hom}_{A^{e}}(K'_{l-r}(A),A^{e})$ and $L'_{r-l}(A)$ are isomorphic. Returning to the
assumptions of the beginning of the section, we have obtained:
$$\mathrm{Ext}_{A^{e}}^{i}(A, A^{e})=H^{i}(L'_{r-l}(A)),\ i\geq 0.$$
Since $A$ is AS-Gorenstein, $H^{i}(L'_{r}(A))=0$ if $i\neq D$. But $L'_{r}(A)= k\otimes_{A}
L'_{r-l}(A)$. Therefore, Proposition 4.1 shows that $\mathrm{Ext}_{A^{e}}^{i}(A, A^{e})=0$ if $i\neq
D$.

Now recall that the Frobenius algebra $E(A)$ (Corollary 5.12) is defined by the $E(A)$-grMod
isomorphism $\varphi: E(A) \rightarrow E(A)^{\ast}(-D)$ such that $\varphi(1)=u^{\ast}$.
The associated Frobenius pairing $(-,-)$~\cite{smi:some} is defined by 
$$(x,y)= \langle x,\ \varphi (y) \rangle = \langle x \bullet y,\ u^{\ast} \rangle, \ x \in E(A), \
y \in E(A).$$ 
For $x \in E(A)_{i}$ and $y\in E(A)_{j}$, $(x,y)=0$ if $i+j\neq D$. Assume $i+j=D$. Then
Proposition 3.1 shows that 
$$(x,y)=(-1)^{i.(D-i)}\langle xy, \ u^{\ast}\rangle,$$
where $xy$ is the product in $A^{!}$. Let us be more explicit in this case. Write down
$u=\bar{\mathtt{u}}$, $\mathtt{u}\in (V^{\otimes n(D)})^{\ast}$. The canonical isomorphism
$A_{n(D)}^{!} \cong W_{n(D)}^{\ast}$ (see end of Section 3) maps $u$ to $\mathtt{u}|_{W_{n(D)}}$.
By the transposed isomorphism, $u^{\ast}$ is the image of an element $w\in W_{n(D)}\subseteq
V^{\otimes n(D)}$. Considering
$w$ as a linear form on $(V^{\otimes n(D)})^{\ast}$, $w$ is the unique element of $W_{n(D)}$ which
maps
$\mathtt{u}$ to 1. Then if $x= \bar{f}$, $y=\bar{g}$, $f\in (V^{\otimes n(i)})^{\ast}$, $g\in
(V^{\otimes n(D-i)})^{\ast}$, we have
\begin{equation}\label{pairing}
(x,y)=(-1)^{i.(D-i)}\langle f\otimes g, \ w\rangle,
\end{equation}
where $f\otimes g$ is the usual tensor product of the linear forms $f$
and $g$. If $N>2$, the sign in (\ref{pairing}) is $+1$.

Recall that the graded Frobenius algebra $E(A)$ is said to be \emph{symmetric} (resp. \emph{graded
symmetric}) if for any $x \in E(A)_{i}$ and $y\in E(A)_{D-i}$, $(y,x)=(x,y)$ (resp.
$(y,x)=(-1)^{i.(D-i)}(x,y)$). With the above notations, it is equivalent to
saying that $\langle g\otimes f, w \rangle= \langle f\otimes g, w \rangle$ (resp.
$\langle g\otimes f, w \rangle=(-1)^{i.(D+1)}\langle f\otimes g, w \rangle$) for any $f$ and $g$. If
$N>2$, both conditions are equivalent.
In the examples, it will be important to compute $w$.

Let us define the automorphism $\nu$ of the graded algebra $E(A)$ by
$$(x,y)=(y,\nu(x)), \ x\in E(A), \ y\in E(A).$$
In particular, one has
\begin{equation}\label{autom}
x\bullet y = y\bullet \nu(x), \ x\in E(A)_{i}, \ y\in E(A)_{D-i}.
\end{equation}  
Clearly $\varphi$ is an isomorphism from the bimodule $E(A)_{\nu}$ to the
bimodule $E(A)^{\ast}$. The notation $E(A)_{\nu}$ means that the right action is twisted by $\nu$
(and the left action is the usual one). In other words, for $a$, $b$,
$x$ in $E(A)$, the action on $x$ in
this bimodule is $(a,b)\mapsto a\bullet x\bullet \nu(b)$. Let $\nu_{1}: V^{\ast}\rightarrow V^{\ast}$ be the
component of degree 1 of $\nu$.
\Blm
$\nu_{1}^{\otimes N}(R^{\perp})$ is included into $R^{\perp}$.
\Elm
\Bdm
From (\ref{autom}), we get the formula
$$e^{\ast}_{i_{1}} \cdots e^{\ast}_{i_{n(D)}}= \pm \, e^{\ast}_{i_{N+1}} \cdots
e^{\ast}_{i_{n(D)}} \nu_{1}(e^{\ast}_{i_{1}}) \cdots
\nu_{1}(e^{\ast}_{i_{N}}).$$
Let $f= \sum \lambda_{i_{1},\ldots , i_{N}} e^{\ast}_{i_{1}}\otimes \cdots \otimes
e^{\ast}_{i_{N}}$ be any element of $R^{\perp}$. The previous formula shows that 
$$x\bullet \overline{\nu_{1}^{\otimes N}(f)}=0$$
for any $x$ in $E(A)_{D-2}$. Thus $\overline{\nu_{1}^{\otimes N}(f)}=0$ since the Frobenius pairing
is non-degenerate. 
\Edm

Define now $\phi_{1}: V \rightarrow V$ as being the transposed linear map of $\nu_{1}$. Lemma
6.1 implies that $\phi_{1}^{\otimes N}(R)$ is included into $R$. So Tens$(\phi_{1})$ defines an
automorphism $\phi$ of the graded algebra $A$, homogeneous of degree 0. Our aim is to compute the
bimodule $U=H^{D}(L'_{r-l}(A))$, i.e., the cokernel of
\begin{equation}
\partial = \partial_{r}+(-1)^{D}\partial_{l} : A\otimes A_{n(D)-1}^{!} \otimes A
\longrightarrow A\otimes A_{n(D)}^{!} \otimes A.
\end{equation}  
Let $\varepsilon$ be the automorphism of $A$ which is the multiplication by
$(-1)^{m}$ on $A_{m}$. Define $\mu_{u}^{\phi}:A\otimes A_{n(D)}^{!} \otimes A \rightarrow A$ by
$\mu_{u}^{\phi}(a\otimes u
\otimes b)=(\varepsilon^{D+1}\phi)(a)b$. Adapting the computation made in the proof of Theorem
9.2 of~\cite{vdb:exist}, we prove the following.
\Blm
$\mu_{u}^{\phi} \circ \partial=0.$
\Elm
\Bdm 
Introduce $(\zeta_{i}^{\ast})_{i}$ the left dual basis to $(\zeta_{i}=e^{\ast}_{i})_{i}$ for the
Frobenius pairing. In other words, $\zeta_{i}^{\ast} \in E(A)_{D-1}$ with
$(\zeta_{i}^{\ast},e^{\ast}_{j})=\delta_{ij}$. One has 
$$\partial(1\otimes \zeta_{i}^{\ast} \otimes 1)= \sum_{j} 1\otimes e^{\ast}_{j}\zeta_{i}^{\ast}
\otimes e_{j} +(-1)^{D} e_{j} \otimes \zeta_{i}^{\ast}e^{\ast}_{j} \otimes 1.$$
Proposition 3.1 gives $\zeta_{i}^{\ast}e^{\ast}_{j}=(-1)^{D-1}
\zeta_{i}^{\ast}\bullet e^{\ast}_{j}= (-1)^{D-1} \delta_{ij} u$. Moreover
$$e^{\ast}_{j}\zeta_{i}^{\ast}= (-1)^{D-1} e^{\ast}_{j}\bullet \zeta_{i}^{\ast}=
(-1)^{D-1}\zeta_{i}^{\ast} \bullet \nu(e^{\ast}_{j})= (-1)^{D-1}(\zeta_{i}^{\ast},\,
\nu(e^{\ast}_{j}))u.$$
But the linear forms $x\mapsto (\zeta_{i}^{\ast},x)$ and $x\mapsto \langle e_{i},x \rangle$ on
$V^{\ast}$ are the same. Therefore $e^{\ast}_{j}\zeta_{i}^{\ast}= (-1)^{D-1} \langle \phi
(e_{i}),e_{j}^{\ast} \rangle u$. So we get
$$\partial(1\otimes \zeta_{i}^{\ast} \otimes 1)=(-1)^{D-1}[ 1\otimes u \otimes \phi(e_{i}) +
(-1)^{D}e_{i} \otimes u \otimes 1].$$
We conclude that
$$\mu_{u}^{\phi} \circ \partial (1\otimes \zeta_{i}^{\ast} \otimes 1)=(-1)^{D-1}[\phi(e_{i})+
(-1)^{D} (-1)^{D+1} \phi(e_{i})]=0.\ \Edm $$
Applying the functor $k\otimes_{A}-$ to the $A$-grMod-$A$ complex
\begin{equation} \label{cok}
A\otimes A_{n(D)-1}^{!} \otimes A \stackrel{\partial}\longrightarrow A\otimes A_{n(D)}^{!} \otimes
A \stackrel{\mu_{u}^{\phi}}\longrightarrow _{\varepsilon^{D+1}\phi}A(n(D)) \rightarrow 0,
\end{equation}  
we get the exact grMod-$A$ complex
$$A_{n(D)-1}^{!} \otimes A \stackrel{\xi_{r}\cdot}\longrightarrow A_{n(D)}^{!} \otimes
A \stackrel{\epsilon_{u}}\longrightarrow k(n(D)) \rightarrow 0.$$
Therefore, Proposition 4.1 shows that (\ref{cok}) is exact, and we have
$$U=_{\varepsilon^{D+1}\phi}A(n(D)).$$
Thus the bimodule $U$ is invertible and $U\otimes_{A}M \cong _{\varepsilon^{D+1}\phi}M$. So we
obtain the following $N$-generalization of Van den Bergh's Proposition 2~\cite{vdb:dual}.
\Bte \label{duality}
Assume that $A$ is an $N$-homogeneous algebra which is generalized Koszul, AS-Gorenstein and with
finite global dimension $D$. Let $\varepsilon$ be the automorphism of $A$ which is the multiplication
by
$(-1)^{m}$ on $A_{m}$. Let $\nu$ be the automorphism of the Frobenius algebra
$E(A)$ such that the bimodule $E(A)^{\ast}$ is canonically $E(A)_{\nu}$. Let $\phi$ be the
automorphism of $A$ such that the homogeneous component $\phi_{1}$ is the transposed linear map
of $\nu_{1}$. Then for any $A$-$A$-bimodule $M$, we have 
$$HH^{i}(A,M)\cong HH_{D-i}(A,\,_{\varepsilon^{D+1}\phi}M).$$
\Ete

Let us examine the condition for having no twist on $M$, i.e.,
$\phi=\varepsilon^{D+1}$. This condition is related to the
existence of a bimodule version of the left $A^{!}\otimes A$-linear map $\Phi : A^{!}\otimes A \rightarrow
A^{!\ast}\otimes A$ defined just before Theorem 5.4. Setting
$\Phi=\Phi_{r}$, recall that it is defined by $\Phi_{r}(1\otimes 1)= u^{\ast}\otimes 1$. Introduce the right $A\otimes A^{!}$-linear map $\Phi_{l} : A\otimes A^{!} \rightarrow
A\otimes A^{!\ast}$ by $\Phi_{l}(1\otimes 1)=1 \otimes u^{\ast}$. Then
$1_{A}\otimes \Phi_{r}$ and $\Phi_{l} \otimes 1_{A}$ are $A$-$A$-bimodule
morphisms from $A\otimes A^{!}\otimes A$ to $A\otimes
A^{!\ast}\otimes A$. Let $\lambda$ be the graded bimodule automorphism of
$A\otimes A^{!}\otimes A$ which is the multiplication by
$(-1)^{m.(D+1)}$ on $A^{!}_{m}$. Notice that $\lambda$ is the identity
automorphism when $N>2$. 
\Bpo
Keep the assumptions of Theorem \ref{duality}. Then
$\phi=\varepsilon^{D+1}$ if and only if 
\begin{equation} \label{morbi}
1_{A}\otimes \Phi_{r}=(\Phi_{l} \otimes 1_{A})\circ \lambda.
\end{equation}
\Epo
\Bdm
(\ref{morbi}) is equivalent to
$f.u^{\ast}=(-1)^{m.(D+1)}u^{\ast}.f$ for any $f$ in $A^{!}_{m}$,
i.e., $e_{i}^{\ast}.u^{\ast}=(-1)^{D+1}u^{\ast}.e_{i}^{\ast}$ for
any $i$. Recall that for $f\in E(A)_{i}$ and $\alpha \in
E(A)^{\ast}_{i+j}$ with $n(i+j)=n(i)+n(j)$, $f\bullet \alpha =(-1)^{ij} f\cdot \alpha$ and $\alpha \bullet f =(-1)^{ij} \alpha \cdot f$. Using (\ref{autom}), it is easy to check that
$\alpha \bullet f= \nu(f)\bullet \alpha$. So (\ref{morbi}) is equivalent to
$$e_{i}^{\ast}\bullet
  u^{\ast}=(-1)^{D+1} \nu (e_{i}^{\ast})\bullet  u^{\ast}$$
for any $i$, that is to $\nu_{1}=(-1)^{D+1} 1_{V^{\ast}}$. 
\Edm
\\

Assume that $\phi=\varepsilon^{D+1}$. Then $\Phi_{r-l}=1_{A}\otimes
\Phi_{r}=(\Phi_{l} \otimes 1_{A})\circ \lambda$ is an $N$-complex
morphism from $L_{r-l}(A)^{ch}[-n(D)]$ to $K_{r-l}(A)$. After the adequate
contraction, we get a 2-complex isomorphism $\Phi'_{r-l}$ from
$L'_{r-l}(A)^{ch}[-D]$ to $K'_{r-l}(A)$. It is an
isomorphim of \emph{resolutions} from
$L'_{r-l}(A)^{ch}[-D] \stackrel{\mu_{u}^{\phi}} \longrightarrow A$ to
$K'_{r-l}(A) \stackrel{\mu} \longrightarrow A$, since
$\mu_{u}^{\phi}(1\otimes u \otimes 1)=1=\mu(1\otimes 1)$. That
isomorphism has been already obtained by the second author for the generic cubic
AS-regular algebras of global dimension three and of type
A, allowing him to prove the Poincar\'e duality for
these algebras when $M=A$~\cite{nic:cubic}.

Now let us give another proof of Theorem \ref{duality}, without using
the Van den
Bergh duality theorem, but for special bimodules $M$. Let $E$ be a
$k$-vector space. Then $M= E\otimes A$ (tensor product over $k$) is naturally a
bimodule. Set $L'_{r-l}(M)=E\otimes L'_{r-l}(A)$. Applying the exact
functor $E\otimes -$ to the resolution
$$L'_{r-l}(A)\stackrel{\mu_{u}^{\phi}}\longrightarrow 
_{\varepsilon^{D+1}\phi}A\rightarrow 0,$$ 
we see that the chain complex
$L'_{r-l}(M)^{ch}[-D]$ is a resolution of the bimodule $_{\varepsilon^{D+1}\phi}M$. Therefore
$$HH_{D-i}(A,\,_{\varepsilon^{D+1}\phi}M)\cong H_{D-i}(L'_{r-l}(M)^{ch}[-D]\otimes_{A^{e}}A)=
H^{i}(L'_{r-l}(M)\otimes_{A^{e}}A).$$
For any $A$-$A$-bimodules $N$ and $P$, the canonical isomorphisms
$$(E\otimes N)\otimes_{A^{e}}P\cong E\otimes (N\otimes_{A^{e}}P)\cong
(N\otimes_{A^{e}}P) \otimes E \cong N\otimes_{A^{e}}(P \otimes E)$$
show that the complex $L'_{r-l}(M)\otimes_{A^{e}}A$ is isomorphic to
$L'_{r-l}(A)\otimes_{A^{e}}M$. Thus
$$H^{i}(L'_{r-l}(M)\otimes_{A^{e}}A)\cong
H^{i}(\mathrm{Hom}_{A^{e}}(K'_{l-r}(A),A^{e})\otimes_{A^{e}}M)$$
$$\cong H^{i}(\mathrm{Hom}_{A^{e}}(K'_{l-r}(A),M))\cong HH^{i}(A,M).$$

As shown by Van den Bergh in the quadratic case (\cite{vdb:exist}, Corollary 9.3), the 
automorphism $\phi$ is related to the Artin-Schelter matrix $Q$
\cite{as:regular} in the cubic case. Let $A$ be
an AS-regular algebra (with polynomial growth) of global dimension 3,
with cubic relations. Then $A$ is Koszul (Proposition 5.2). The basis
of $V$ is $(x,y)$, and the basis of $R$ is $(f_{1}, f_{2})$. The nonvanishing components of
$E(A)$ are: $E(A)_{0}=k$, $E(A)_{1}=V^{\ast}$ with basis $(x^{\ast},
y^{\ast})$, $E(A)_{2}=R^{\ast}$ with basis $(f_{1}^{\ast},f_{2}^{\ast})$, 
$E(A)_{3}= W_{4}^{\ast}$ with basis $w^{\ast}$. We have
\begin{equation}\label{w}
w=xf_{1} +y f_{2} = g_{1}x + g_{2}y,
\end{equation} 
where $(g_{1}, g_{2})$ is the basis of $R$ defined by the matrix product
$$(g_{1}\ g_{2})= (f_{1}\ f_{2}) Q^{t}.$$
Here $Q=(q_{ij})$ is some invertible $2\times 2$ matrix with scalar
entries. Using (\ref{pairing}) and (\ref{w}), it is easy to compute
the Frobenius pairing. We find
$$(x^{\ast}, f_{1}^{\ast})=(y^{\ast}, f_{2}^{\ast})=1, \ (x^{\ast},
f_{2}^{\ast})=(y^{\ast}, f_{1}^{\ast})=0,$$
$$(f_{1}^{\ast}, x^{\ast})=q_{11},\ (f_{2}^{\ast}, x^{\ast})=q_{12},\
(f_{1}^{\ast}, y^{\ast})=q_{21},\ (f_{2}^{\ast}, y^{\ast})=q_{22}.$$
We draw $\nu^{-1}(x^{\ast}) \bullet f_{1}^{\ast}=
f_{1}^{\ast} \bullet x^{\ast} =q_{11}w^{\ast}$, and next,
$\nu^{-1}(x^{\ast}) \bullet f_{2}^{\ast} =q_{12}w^{\ast}$. Thus
$\nu^{-1}(x^{\ast})= q_{11}x^{\ast} +q_{12} y^{\ast}$. Similarly,
$\nu^{-1}(y^{\ast})= q_{21}x^{\ast} +q_{22} y^{\ast}$. Denoting by
$X^{\ast}$ the column vector with two entries $x^{\ast}$ and
$y^{\ast}$, the isomorphism $\nu_{1}$ has matrix form $X^{\ast}
\mapsto Q^{-1}X^{\ast}$. Denoting by
$F^{\ast}$ the column vector with two entries $f_{1}^{\ast}$ and
$f_{2}^{\ast}$, a similar computation shows that the
isomorphism $\nu_{2}$ has matrix form $F^{\ast} \mapsto
Q^{t} F^{\ast}$. Thus the Frobenius algebra $E(A)$ is symmetric if and
only if $A$ is of type A. We have also obtained the following.
\Bpo
Assume that the algebra $A$ is AS-regular (with polynomial growth), of global dimension $3$,
with cubic relations.
The isomorphism $\phi_{1}$ has matrix form $X
\mapsto (Q^{-1})^t X$. In particular,
$\varepsilon^{D+1}\phi=1_{A}$ if and only if $A$ is of type \emph{A}.
\Epo

We finish by another class of examples. Consider $A$ as in Corollary
5.10 with $N=2$ or ($N>2$ and $n=Nq+1$). Then $D=n$ or $D=2q+1$,
respectively. In both cases, $n(D)=n$. Recall that
$u=x_{n}^{\ast}\cdots x_{1}^{\ast}$. A straightforward computation provides
$$w=\sum_{\sigma \in \mathbf{S}_{n}} \mathrm{sgn} (\sigma)\ x_{\sigma(n)} \otimes \cdots
\otimes x_{\sigma(1)}\,$$
where $\mathbf{S}_{n}$ denotes the permutation group of $\{1, \ldots ,n\}$. We leave to the
reader the computation of the Frobenius pairing by using
(\ref{pairing}). For $j_{1}>\cdots >
j_{n(i)}$, $0\leq i \leq D$, one has
$$\nu (x^{\ast}_{j_{1}} \cdots 
x^{\ast}_{j_{n(i)}})=(-1)^{n(i).(n+1)}x^{\ast}_{j_{1}} \cdots
 x^{\ast}_{j_{n(i)}}.$$ 
Thus $E(A)$ is symmetric if and only if $n$ is odd, and for $N=2$, it is
always graded symmetric. On the other hand,
$\phi=\varepsilon^{n+1}$, so that the equality 
$\varepsilon^{D+1}\phi=1_{A}$ holds for $N=2$, and this equality holds for
$N>2$ if and only if $n$ is odd.

\end{document}